\documentclass{article}
\usepackage{amssymb,amsmath,amsfonts}
\usepackage{color,soul}
\usepackage{graphicx}
\usepackage{latexsym}
\usepackage{subfigure}
\usepackage{relsize}
\usepackage{amsthm}
\usepackage{amsmath, amsthm, amscd, amsfonts, amssymb}
\usepackage{cite}
\newtheorem{thm}{Theorem}[section]
 
 \newtheorem{lem}[thm]{Lemma}

 \newtheorem{defn}[thm]{Definition}
 \makeatletter
\def\@seccntformat#1{\csname the#1\endcsname.\quad}
\def\numberline#1{\hb@xt@\@tempdima{#1\if&#1&\else.\fi\hfil}}
\makeatother

\begin{document}
\title{\textbf{Fractional Order Runge-Kutta Methods}}
\author{\small{\bf R. Ghaffari$^{1\star}$,  F. Ghoreishi$^{2\star}$ and Nasser {Saad$^{3\star}${\footnote{Corresponding
Author}}} }\\
 \footnotesize{$1$ Faculty of Mathematics, K. N. Toosi University of Technology, Tehran, Iran.; r.ghaffari85@gmail.com} \\
 \footnotesize{$2$ Faculty of Mathematics, K. N. Toosi University of Technology, Tehran, Iran.; ghoreishif @kntu.ac.ir} \\
 \footnotesize{$3$ School of Mathematical and Computational Sciences, University of Prince Edward Island,}\\
\footnotesize{Charlottetown,Canada.; nsaad@upei.ca}\\
\footnotesize{$\star$ These authors contributed equally to this work.}
}\maketitle
\date{}
\maketitle%
%%%%%                         ABSTRACT
%%%%%%%%%%%%%%%%%%%%%%%%%%%%%%%%%%%%%%%%%%%%%%%%%%%%%%%%% %
\begin{abstract}
 {This paper investigates, a new class of fractional order Runge-Kutta
(FORK) methods for numerical approximation to the solution of fractional
differential equations (FDEs). By using the Caputo generalized Taylor
formula for Caputo fractional derivative, we construct explicit and implicit
FORK methods, as the well-known Runge-Kutta schemes for ordinary
differential equations. In the proposed method, due to the dependence of
fractional derivatives to a fixed base point $t_0$, we had to modify the right-hand side of the given equation in all steps of the FORK methods. Some
coefficients for explicit and implicit FORK schemes are presented. The
convergence analysis of the proposed method is also discussed. Numerical
experiments clarify the effectiveness and robustness of the method..}
\end{abstract}
 $\bf{Keywords }$ Fractional differential equations; Caputo fractional derivative; Convergence analysis; Consistency;
 Stability analysis.


 %%%%%%%%%%%%%%%%% INTRODUCTION%%%%%%%%%%%%%%%%%%%%%%%%%%%%%%%%
 %%%%%%%%%%%%%%%%%%%%%%%%%%%%%%%%%%%%%%%%%%%%%%%%%%%%%%%%%%%
 \section{Introduction}

In recent years, {the numerical approximation for the solutions of} FDEs has attracted increasing attention in many fields of applied sciences and engineering \cite{1,7,8}.  It is common for FDEs to be used in formulating many problems in applied mathematics. Developing numerical methods for fractional differential problems is necessary and important because analytic solutions are usually challenging to obtain. Moreover, it is necessary to develop numerical methods that are highly accurate and easy to use.

It is well known that fractional derivatives have different definitions; the most common and important ones in applications are the Riemann--Liouville and Caputo fractional derivatives. Models describing physical phenomena usually prefer the use of the Caputo derivative. One of the reasons is that the Riemann--Liouville derivative needs initial conditions containing the limit values of the Rieman--Liouville fractional derivative at the origin of time. In contrast, the initial conditions for Caputo derivatives are the same as for integer-order differential equations. Therefore, using the Caputo derivative, there is a clear physical interpretation of the prescribed data; see \cite{1,17,18}.

Numerous research papers have been published on numerical methods for FDEs. Many researchers considered the trapezoidal method, predictor-corrector method, extrapolation method, and spectral method \cite{41,42,43,44,45,451,452,9,10,11,191}. Some of these methods discretize fractional derivatives directly. As an example, the L1 formula was created by a piecewise linear interpolation approximation for the integrand function on each small interval \cite{19,20}. In \cite{16}, the authors applied quadratic interpolation approximation using three points to approximate the Caputo fractional derivative, while in \cite{15}, a technique based on the block-by-block approach was presented. This technique became a common method for equations with integral operators. In \cite{21}, Caputo fractional differentiation
was approximated by a weighted sum of the integer order derivatives of functions. In \cite{14}, several numerical algorithms were proposed to approximate the Caputo fractional derivatives by applying higher-order piecewise interpolation polynomials and the Simpson method to design a higher-order algorithm.

{These methods are appropriate options if the resulting system of equations, generated from the numerical method, is linear and well-conditioned. However, they present a high computational cost when the problem we are solving is badly conditioned or nonlinear.
In light of the above discussion and the analysis of other methods for FDEs, despite many papers on numerical methods for FDEs, there are still insufficient efficient numerical approaches for such equations. Therefore, further studies are still in demand. In this case, step-by-step methods such as the Runge--Kutta method are a good option. They are favored due to their simplicity in both calculation and analysis. }

Several authors have used Runge--Kutta methods to solve ordinary, partial differential, and integral equations \cite{22, 23,24, 25,26,27,28,29}.
Lubich and others have done some fundamental works regarding Runge--Kutta methods for Volterra integral equations \cite{27,28,29}. They used the order conditions to derive various Runge--Kutta methods.

{One of the efficient implicit Runge--Kutta methods for the numerical approximation of some linear partial differential equations is the Rosenbrock procedure. It is a class of semi-implicit Runge--Kutta methods for the numerical solution of some stiff systems of ODEs. In Osterman and Rochet's papers \cite{241,242}, the authors apply the Rosenbrock methods to solve linear partial differential equations, obtaining a sharp lower bound for the order of convergence. They show that the order of convergence is, in general, fractional. So, for the numerical solution of some fractional linear
partial differential equations, we can construct fractional Rosenbrock-type methods, in which a special type of fractional semi-implicit Runge--Kutta method could be considered. }

{
This paper introduces a new class of fractional order Runge--Kutta methods for numerical approximation to the solution of FDEs. Using the Caputo generalized Taylor series formula for the Caputo fractional derivative, we construct explicit and implicit FORK methods comparable to the well-known Runge--Kutta schemes for ordinary differential equations.}

The remainder of the paper is organized as follows. In Section \ref{P}, we review some definitions and properties of fractional calculus.  We propose new explicit and implicit FORK methods for solving FDEs in Sections \ref{se1} and \ref{se2}.  In Section \ref{st}, the theoretical analysis of the  convergency, stability, and consistency of the proposed methods is presented. Finally, in Section \ref{exam}, some numerical examples demonstrate the effectiveness of the methods proposed. Also, in Appendix A. two Mathematica computer programming codes are given.

  %%%%%%%%%%%%%%%%%%%%%%%%%%%%%%%%%%%%%
\section{Preliminaries}\label{P}
In this section, we briefly state definitions of fractional integral
and Caputo fractional derivative and some of them properties. For
further information about fractional calculus and some other
 definitions of fractional derivatives, we refer the interested readers to the \cite{1,17,188}.
\begin{defn}
The Riemann Liouville fractional integral operator of order $\alpha
>0$ for a function $f(x)\in L_1[a,b]$ with $a\geq 0$ is defined as
\begin{equation*}
J^{\alpha}_af(x)=\frac{1}{\Gamma(\alpha)}\int^{x}_{a}(x-t)^{\alpha
-1}f(t)dt, \quad x\in[a,b], \qquad J^0_af(x)=f(x).
\end{equation*}
where  $L_1[a,b]=\{f~ |~  f ~ \text{is a measurable function on}
 ~[a,b]~ \text{and} ~ \int_a^b |f(x)| dx<\infty\},$  $\Gamma$ is Gamma
function and $\Gamma(\alpha+1)=\alpha !$ \,.
\end{defn}
 %%%%%%%%%%%%%%%%%%%
\begin{defn}\label{d1}
The Caputo fractional derivatives of order $ \alpha >0$ of a
function $f(x)\in L_1[a,b]$ with $a\geq 0$ is defined as
\begin{equation}\label{E2}
{}_a^cD_t^\alpha f(x) =J_a^{n-\alpha}D^nf(x)= \left\{ \begin{gathered}
  \frac{1}{{\Gamma (n - \alpha )}}\int_a^x {{(x - t)}^{n - \alpha  - 1}}D^nf(t)dt,
  ~ n-1<\alpha< n,~ n \in \mathbb{N},\\
  D^nf(x),\begin{array}{*{20}{c}}
  {}&{\alpha  = n.}
\end{array} \hfill \\
\end{gathered}  \right.
\end{equation}
\end{defn}
%%%%%%%%%%%%%%%%%%%%%%%%%%%%%%%%%%%%%%%%%%%%%%%%%%%%%%%%%%%
   \begin{thm}
 (Generalized Taylor formula for Caputo fractional derivative \cite{2}). Suppose that $({}_a^cD_t^{\alpha})^k f(x) \in C(a,b] $ for
 $k=0,1,...,n+1$, where $0<\alpha \leqslant 1$,
 then for $\forall x \in [a,b]$ there exist $\xi \in (a,x)$ such
 that
 \begin{equation}\label{E3}
 f(x)=\sum_{i=0}^{n}\frac{(x-a)^{i\alpha}}{\Gamma(i\alpha +1)}(({}_a^cD_t^{\alpha})^i f)(a)+
  \frac{(x-a)^{(n+1)\alpha}}{\Gamma((n+1)\alpha +1)}(({}_a^cD_t^{\alpha})^{(n+1)}) f)(\xi),
 \end{equation}
 where  $({}_a^cD_t^{\alpha})^n={}_a^cD_t^{\alpha}\, {}_a^cD_t^{\alpha} ...{}_a^cD_t^{\alpha}
 $(n times).
 \end{thm}
There are also two important functions in fractional calculus.
They are direct generalization of the exponential series which play important roles in the solution of FDEs and stability analysis.
\begin{defn}
The Mittag-Leffler function is defined as
\begin{equation*}\label{E4}
E_{\alpha}(x)=\sum_{k=0}^{\infty}\frac{x^k}{\Gamma(\alpha k+1)}, \quad \Re (\alpha)>0, x\in \mathbb{C}.
\end{equation*}
Also, the two-parameters Mittag-Leffler function is defined by
\begin{equation*}\label{E5}
E_{\alpha, \beta}(x)=\sum_{k=0}^{\infty}\frac{x^k}{\Gamma(\alpha k+\beta)}, \quad \Re (\alpha)>0, \beta \in \mathbb{C}, x \in \mathbb{C},
\end{equation*}
\end{defn}
We note that $E_{\alpha}(x)=E_{\alpha, 1}(x)$ and
$$E_1(x)=\sum_{k=0}^{\infty}\frac{x^k}{\Gamma(k+1)}=\sum_{k=0}^{\infty}\frac{x^k}{k!}=exp(x).$$

%%%%%%%%%%%%%%%%%%%%                 Modified Runge Kutta methods %%%%%%%%%%%%%%%%%%%%%%%%%%%%%%%%%%%%%%%%%%%%%%%%%%%%%%%%%%%%%%%%%%%%%%%%%%%%%%%

\section{Fractional order Runge-Kutta methods}\label{se1}
In this section, a new class of FORK methods for numerical solution
of FDE is investigated.\\

 Consider the following FDE with $0< \alpha
\leq1$:
\begin{eqnarray}
&& {}_{t_0}^cD_t^{\alpha} y(t)=f(t,y(t)),\quad t\in[ t_0,T],\nonumber\\
&& y(t_0)=y_0.\label{E6}
\end{eqnarray}
where $y(t)\in C[t_0,T]$ and  $f(t,y(t))\in C[t_0,T]\times
\textsl{R}.$ ~  $t_0$ is called base point of fractional derivative.\\
We set $t_n=t_0+nh$, $n=0,1,...,N^m$, where $h= \frac{T-t_0}{N^m} $ is the
step size, $N$ is a positive integer and in section \ref{st}, we will prove that $m\geq\frac{1}{\alpha}$.

For the existence and uniqueness of solution of the FDE (\ref{E6}), consider the following theorem from
\cite{17}.

\begin{thm}\label{Tu}
 Let $\alpha>0$, $y_0\in \mathbb{R}$, $K>0$ and $T>0$ and also let
the function $f:G\rightarrow \mathbb{R}$ be continuous and fulfil a
Lipschitz condition with respect to the second variable, i.e.
$$ |f(t,y_1)-f(t,y_2)|\leq L|y_1-y_2|$$
with some constant $L>0$ independent of $ t, y_1$ and $y_2$.
Define\\
$G=\left\lbrace (t,y):t\in [0,T], |y-y_0|\leq K  \right\rbrace $,
$M=Sup_{(t,z)\in G}|f(t,z)|$ and
\[T^* = \left\{ \begin{gathered}
  \begin{array}{*{20}{c}}
  {T,}&{~~~~~~~~~~~~~~}&{~~~~~~~~~~~~~~~}&{M = 0,}
\end{array} \hfill \\
  \begin{array}{*{20}{c}}
  {\min \{ T,{{(K\Gamma (\alpha + 1)/M)}^{1/\alpha}}\}, }&{}&{else.}&{}
\end{array} \hfill \\
\end{gathered}  \right.\]
Then, there exists a uniquely function  $y\in C[0,T^*]$solving the
initial-value problem (\ref{E6}).
\end{thm}
In the sequel, we assume that  $f(t,y)$ has continuous partial
derivatives with respect to $t$ and $y$ of as high an order as we
desire.

Now, we introduce a s-stage explicit fractional order Runge-Kutta (EFORK) method for FDEs, which is discussed completely with 2 and 3 stages.
%%%%%%%%%%%%%%%%%
\begin{defn}
A family of s-stage EFORK methods is defined as
\begin{eqnarray}
&& K_1=h^{\alpha}f(t,y),\nonumber\\
&& K_2=h^{\alpha}f(t+c_2h,y+a_{21}K_1),\nonumber\\
&& K_3=h^{\alpha}f(t+c_3h,y+a_{31}K_1+a_{32}K_2),\nonumber\\
&&\qquad\vdots \nonumber\\
&& K_s=h^{\alpha}f(t+c_sh,y+a_{s1}K_1+a_{s2}K_2+\cdots +a_{s,s-1}K_{s-1}),\label{E8}
\end{eqnarray}
with
\begin{equation}\label{E9}
y_{n+1}=y_n+\sum_{i=1}^{s}w_iK_i,
\end{equation}
where the unknown coefficients $\left\lbrace a_{ij}\right\rbrace
_{i=2,j=1}^{s,i-1}$ and
 the unknown weights $\left\lbrace c_i\right\rbrace_{i=2}^s $, $\left\lbrace w_i \right\rbrace_{i=1}^s
 $ have to be determined.
\end{defn}

To specify a particular method, one needs to provide $\left\lbrace
a_{ij}\right\rbrace _{i=2,j=1}^{s,i-1}$ and
 $\left\lbrace c_i\right\rbrace_{i=2}^s $, $\left\lbrace w_i \right\rbrace_{i=1}^s $ accordingly. Following Butcher \cite{22},
 a method of this type is designated by the following scheme.\\
 \begin{table}[!h]
\begin{center}
\begin{tabular}{c | c c c c c}
$c_2$      &\quad$a_{21}$\\
$c_3$     &\quad$a_{31}$   &\quad$a_{32}$\\
$\vdots$    &\quad$\vdots$ &\quad$$  &\quad$\ddots$\\
$c_s$     &\quad$a_{s1}$   &\quad$a_{s2}$      &\quad $\cdots$       &\quad$a_{s,s-1}$\\
\hline &\quad $w_1$ &\quad $w_2$     &\quad $\cdots$   &\quad $w_{s-1}$&\quad $w_s$ \\
 \end{tabular}
\end{center}
\end{table}
We expand $y_{n+1}$ in (\ref{E9}), in powers of $h^{\alpha}$, such
that it agrees with the Taylor series expansion of the solution of
the FDE (\ref{E6}) in a specified number of terms, (see \cite{4}).
To do this,
we need some changes on Caputo Taylor series expansion of $y(t)$.\\
According to (\ref{E3}), generalized Taylor formula for $\alpha \in
(0,1]$ with respect to
 Caputo fractional derivative of function $y(t)$ is defined as
follows
\begin{equation}\label{E12}
y(t)=y(t_0)+\frac{(t-t_0)^\alpha }{\Gamma(\alpha+1)}{}_{t_0}^cD_t^{\alpha}y(t_0)+\frac{(t-t_0)^{2\alpha} }{\Gamma(2\alpha+1)}(({}_{t_0}^cD_t^{\alpha})^2y)(t_0)
+\frac{(t-t_0)^{3\alpha} }{\Gamma(3\alpha+1)}(({}_{t_0}^cD_t^{\alpha})^3y)(t_0)+\cdots
\end{equation}
where
\begin{align}
&{}_{t_0}^cD_t^{\alpha}y(t)=f(t,y),\nonumber\\
&({}_{t_0}^cD_t^{\alpha})^2y(t)={}_{t_0}^cD_t^{\alpha}f(t,y),\nonumber\\
&({}_{t_0}^cD_t^{\alpha})^3y(t)=({}_{t_0}^cD_t^{\alpha})^2f(t,y),\nonumber\\
&\qquad \vdots \label{dd1}
\end{align}
For obtaining an explicit expression for (\ref{dd1}), we propose
total differential theorem for Caputo fractional derivative by using
(\ref{E2}) and total differential theorem for derivatives of integer
order.

%%%%%%%%%%%%%%%%%%%
Caputo fractional derivatives of composite function $f(t,y(t))$ can
be computed by fractional Taylor series:
\begin{align}\label{E11}
  f(t,y(t))=f(t_0,y(t_0))+\frac{(t-t_0)^\alpha }{\Gamma(\alpha+1)}f^{\alpha}_{t}(t_0,y(t_0))+\frac{(y-y_0)}{1!}f_y(t_0,y(t_0))\nonumber\\
 + \frac{(t-t_0)^{2\alpha} }{\Gamma(2\alpha+1)}f^{\alpha,\alpha}_{t,t}(t_0,y(t_0)) +\frac{(y-y_0)^{2}}{2!}f_{y,y}(t_0,y(t_0))\hspace{2cm}\nonumber\\
 +\frac{(t-t_0)^\alpha(y-y_0) }{\Gamma(\alpha+1)}f^{\alpha,1}_{t,y}(t_0,y(t_0))+\frac{(t-t_0)^{3\alpha}
}{\Gamma(3\alpha+1)}f^{\alpha,\alpha,\alpha}_{t,t,t}(t_0,y(t_0))+\cdots
\end{align}
where  $f^{\alpha}_t$ is Caputo fractional derivative of $f(t,y(t))$
with respect to $t$. After inserting $y(t)-y(t_0)$ from (\ref{E12}) in (\ref{E11}) and by using fractional derivative of (\ref{E11}) for $\alpha\in (0,1]$, we have
\begin{align*}
{}_{t_0}^cD^{\alpha}_t
f(t,y(t))=f^{\alpha}_{t}(t_0,y(t_0))+f(t,y(t))f_y(t_0,y(t_0))+\frac{(t-t_0)^{\alpha}
}{\Gamma(\alpha+1)}f^{\alpha,\alpha}_{t,t}(t_0,y(t_0))\hspace{-1.3 cm}\nonumber\\
+\frac{1}{2}\left(\frac{\Gamma(2\alpha+1)}{\Gamma(\alpha+1)^3}(t-t_0)^\alpha f^2(t_0,y(t_0))+\ldots\right)f_{y,y}(t_0,y(t_0))\hspace{-.31 cm}\nonumber\\+
\left(\frac{\Gamma(2\alpha+1)}{\Gamma(\alpha+1)^3}(t-t_0)^\alpha f(t_0,y(t_0))+\ldots\right)f^{\alpha,1}_{t,y}(t_0,y(t_0))
\nonumber\\+
\frac{(t-t_0)^{2\alpha}
}{\Gamma(2\alpha+1)}f^{\alpha,\alpha,\alpha}_{t,t,t}(t_0,y(t_0))+\cdots,\hspace{3.15 cm}
\end{align*}
and so
\begin{align}\label{fd1}
{}_{t_0}^cD^{\alpha}_t
f(t_0,y(t_0))=f^{\alpha}_{t}(t_0,y(t_0))+f(t_0,y(t_0))f_y(t_0,y(t_0)).
\end{align}
Also
\begin{align*}
\hspace{-2 cm}{}_{t_0}^cD^{2\alpha}_t f(t,y(t))={}_{t_0}^cD^{\alpha}_t
f(t_0,y(t_0)) f_y(t_0,y(t_0))+
f^{\alpha,\alpha}_{t,t}(t_0,y(t_0))
\nonumber\\+\left(\frac{\Gamma(2\alpha+1)}{2\Gamma(\alpha+1)^2}f^2(t_0,y(t_0))+\ldots\right)f_{y,y}(t_0,y(t_0))\hspace{-.1 cm}\nonumber\\+
\left(\frac{\Gamma(2\alpha+1)}{\Gamma(\alpha+1)^2}f(t_0,y(t_0))+\ldots\right)f^{\alpha,1}_{t,y}(t_0,y(t_0))\hspace{.15 cm}\nonumber\\+\frac{(t-t_0)^{\alpha}
}{\Gamma(\alpha+1)}f^{\alpha,\alpha,\alpha}_{t,t,t}(t_0,y(t_0))+\cdots\hspace{2.2 cm},
\end{align*}
which yields
\begin{align}\label{fd2}
{}_{t_0}^cD^{2\alpha}_t f(t_0,y(t_0))=
f^{\alpha}_{t}(t_0,y(t_0))f_y(t_0,y(t_0))+f(t_0,y(t_0))f^2_y(t_0,y(t_0))+
f^{\alpha,\alpha}_{t,t}(t_0,y(t_0))\nonumber\\
+\frac{1}{2}f^2(t_0,y(t_0))f_{y,y}(t_0,y(t_0))+
f(t_0,y(t_0))f^{\alpha,1}_{t,y}(t_0,y(t_0)).\hspace{1.9 cm}
\end{align}

In a similar manner with (\ref{fd1}-\ref{fd2}), We can obtain the
higher fractional derivatives of $f(t,y(t))$.

%%%%%%%%%%%%%%%%%%%%%%%%%%%%%%%%%%%%%%%%%%%%
Now, by using (\ref{fd1}-\ref{fd2}), we have
\begin{align}\label{EE111}
{}_{t_0}^cD_t^{\alpha}y(t_0)=f(t_0,y_0),\hspace{10 cm}\nonumber\\
({}_{t_0}^cD_t^{\alpha})^2y(t_0)=f^{\alpha}_t(t_0,y_0)+f(t_0,y_0)f_y (t_0,y_0)\hspace{6.5 cm},\nonumber\\
({}_{t_0}^cD_t^{\alpha})^3y(t_0)=f^{\alpha}_{t}(t_0,y(t_0))f_y(t_0,y(t_0))+f(t_0,y(t_0))f^2_y(t_0,y(t_0))+
f^{\alpha,\alpha}_{t,t}(t_0,y(t_0))\hspace{1 cm}\nonumber\\+\frac{1}{2}f^2(t_0,y(t_0))f_{y,y}(t_0,y(t_0))+
f(t_0,y(t_0))f^{\alpha,1}_{t,y}(t_0,y(t_0)).\hspace{3 cm}
\end{align}
where $f^{\alpha, i}_{t,y}$, $i=1,2,...$ represents the $i$th
integer derivative of the function $f_t^{\alpha}$ with respect to
$y$. As we can see from (\ref{E12}), in Caputo fractional
derivatives $(({}_{t_0}^cD_t^{\alpha})^k y)(t_0)$, $k=0,1,2,...$,
argument $t_0$ in $y(t_0)$ and starting value in
$({}_{t_0}^cD_t^{\alpha})^k$ are the same. To construct an efficient
numerical scheme, we should obtain a similar series with the
derivatives evaluated in any other point ($t_n>t_0$), such that the
expansion can be constructed independently from
 the starting point $t_0$. In other words, we need
\begin{align}\label{EE13}
y(t_{n+1})&=y(t_n)+\frac{h^\alpha
}{\Gamma(\alpha+1)}{}_{t_n}^cD_t^{\alpha}y(t_n)+\frac{h^{2\alpha}
}{\Gamma(2\alpha+1)}(({}_{t_n}^cD_t^{\alpha})^2y)(t_n)\nonumber\\
&~~~~~~~~~~~+\frac{h^{3\alpha}
}{\Gamma(3\alpha+1)}(({}_{t_n}^cD_t^{\alpha})^3y)(t_n)+\cdots,
\end{align}
and so $(({}_{t_n}^cD_t^{\alpha})^iy)(t_n), i=1,2,...$. To do so, by
using ${}_{t_0}^cD_t^{\alpha}y(t)$, we  obtain
${}_{t_n}^cD_t^{\alpha}y(t)$ for $n=1,2,...,N^m-1$, as
\begin{align}
{}_{t_n}^cD_t^{\alpha}y(t)&={}_{t_0}^cD_t^{\alpha}y(t)-\frac{1}{\Gamma(1-\alpha)}\int_{t_0}^{t_n}(t-s)^{-\alpha}Dy(s)ds \nonumber\\
&={}_{t_0}^cD_t^{\alpha}y(t)-\frac{1}{\Gamma(1-\alpha)}\sum_{i=0}^{n-1}\int_{t_i}^{t_{i+1}}(t-s)^{-\alpha}Dy(s)ds.
\label{EE16}
\end{align}
By using Lagrange interpolation formula for $y(s)$ in support points
$\{t_i,t_{i+1}\}$, we have
\begin{align*}
y(s)&\simeq\frac{(s-t_i)}{(t_{i+1}-t_i)}y_{i+1}-\frac{(s-t_{i+1})}{(t_i-t_{i+1})}y_i\\
&=\frac{(s-t_i)}{h}y_{i+1}-\frac{(s-t_{i+1})}{h}y_i, \quad s\in[t_i,t_{i+1}],\quad i=0,1,...,n-1,
\end{align*}
where for sufficiently small h we have
\begin{equation*}
Dy(s)\simeq \frac{1}{h}(y_{i+1}-y_i), \quad s\in[t_i,t_{i+1}], \quad
i=0,1,...,n-1,
\end{equation*}
and
\begin{equation*}
\int_{t_i}^{t_{i+1}}(t-s)^{-\alpha}Dy(s)ds \simeq
\frac{(y_{i+1}-y_i)}{h(1-\alpha)}\left[
(t-t_i)^{1-\alpha}-(t-t_{i+1})^{1-\alpha}\right], \quad
i=0,1,...,n-1.
\end{equation*}
From (\ref{EE16}) and ${}_{t_0}^cD_t^{\alpha}y(t)=f(t,y)$, we have
\begin{equation*}\label{h1}
{}_{t_n}^cD_t^{\alpha}y(t)=f(t,y)-\frac{1}{\Gamma(1-\alpha)}\sum_{i=0}^{n-1}\frac{(y_{i+1}-y_i)}{h(1-\alpha)}\left[
(t-t_i)^{1-\alpha}-(t-t_{i+1})^{1-\alpha}\right] .
\end{equation*}

So we may write
\begin{equation}\label{h1}
{}_{t_n}^cD_t^{\alpha}y(t)=F_n(t,y),~~~~~~~ n=0,1,2,... .
\end{equation}
where $F_0(t,y)=f(t,y)$ and for $n=1,2,3,...$ we have
$$F_n(t,y)=f(t,y)-\frac{1}{\Gamma(2-\alpha)}\sum_{i=0}^{n-1}\frac{y_{i+1}-y_i}{h}\left[
(t-t_i)^{1-\alpha}-(t-t_{i+1})^{1-\alpha}\right].$$

Clearly $F_n(t,y)$,  is continuous and Lipschitz
condition with respect to the second variable, due to such
properties of $f(t,y).$ In what follows, for convenience of notation
we rename $F_n(t,y)$ as $f(t,y)$, i.e., in any initial points
$t_n>t_0$, $n=1,2,...,N^m-1$, we consider the right terms of
(\ref{h1}) as $f(t_n,y_n)$ instead $F_n(t_n,y_n)$ in any stages.

Now, for constructing FORK methods, we can use the Taylor formula
(\ref{E12}) and (\ref{EE111}), where ${}_{t_n}^cD_t^{\alpha}y(t)$ is
defined in (\ref{h1}).

%%%%%%%%%%%                  2-stage

\subsection{EFORK method of order $2 \alpha$}
Let us introduce following EFORK method with 2-stage:
\begin{eqnarray}
&& K_1=h^{\alpha}f(t_n,y_n),\nonumber\\
&& K_2=h^{\alpha}f(t_n+c_2h,y_n+a_{21}K_1),\nonumber\\
&&y_{n+1}=y_n+w_1K_1+w_2K_2\label{E10}
\end{eqnarray}
where coefficients $c_2, a_{21}$ and weights $w_1, w_2$ are chosen
to make approximate value $y_{n+1}$ as possible as closer to exact
value  $y(t_{n+1})$. We expand $K_1$ and $K_2$ about the point
$(t_n,y_n)$, where we use Caputo Taylor formula  (\ref{EE13}) about
point $t_n$ and standard integer order Taylor formula about $y_n$ as
\begin{eqnarray*}
&&K_1=h^{\alpha}f(t_n,y_n),\\
&&K_2=h^{\alpha}f(t_n+c_2h,y_n+a_{21}K_1)\\
&&\qquad=h^{\alpha}[ f(t_n,y_n)+\frac{c_2^{\alpha}h^{\alpha}}{\Gamma(\alpha+1)}f_t^{\alpha}+
a_{21}h^{\alpha}f_nf_y+\frac{c_2^{2\alpha}h^{2\alpha}}{\Gamma(2\alpha+1)}f_{t,t}^{\alpha,\alpha}+\frac{a_{21}^2h^{2\alpha}}{2}f_n^2f_{y,y}\\
&&\qquad
+\frac{c_2^{\alpha}a_{21}h^{2\alpha}}{\Gamma(\alpha+1)}f_nf_{t,y}^{\alpha,1}+\cdots
]
=h^{\alpha}f_n+h^{2\alpha}\left( \frac{c_2^{\alpha}}{\Gamma(\alpha+1)}f_t^{\alpha}+a_{21}f_nf_y\right) \\
&&\qquad+h^{3\alpha}\left(
\frac{c_2^{2\alpha}}{\Gamma(2\alpha+1)}f_{t,t}^{\alpha,\alpha}+\frac{a_{21}^2}{2}f_n^2f_{y,y}
+\frac{c_2^{\alpha}a_{21}}{\Gamma(\alpha+1)}f_nf_{t,y}^{\alpha,1}\right) +\cdots
\end{eqnarray*}
Substituting $K_1$ and $K_2$ in (\ref{E10}), we have
\begin{align}
y_{n+1}&=y_n+(w_1+w_2)h^{\alpha}f_n+h^{2\alpha}w_2\left( \frac{c_2^{\alpha}}{\Gamma(\alpha+1)}f_t^{\alpha}+a_{21}f_nf_y\right)\nonumber \\
&\qquad +w_2h^{3\alpha}\left(
\frac{c_2^{2\alpha}}{\Gamma(2\alpha+1)}f_{t,t}^{\alpha,\alpha}+\frac{a_{21}^2}{2}f_n^2f_{y,y}
+\frac{c_2^{\alpha}a_{21}}{\Gamma(\alpha+1)}f_nf_{t,y}^{\alpha,1}\right) +\cdots\label{E13}
\end{align}
Comparing (\ref{EE13}) with (\ref{E13}) and matching coefficients of powers of $h^{\alpha}$, we obtain three equations
\begin{eqnarray}
&&w_1+w_2=\frac{1}{\Gamma(\alpha+1)},\nonumber \\
&&w_2\frac{c_2^{\alpha}}{\Gamma(\alpha+1)}=\frac{1}{\Gamma(2\alpha+1)},\nonumber \\
&&w_2a_{21}=\frac{1}{\Gamma(2\alpha+1)}.\label{EE14}
\end{eqnarray}
From these equations, we see that, if $c_2^{\alpha}$ is chosen arbitrarily (nonzero), then
\begin{equation}\label{E18}
a_{21}=\frac{c_2^{\alpha}}{\Gamma(\alpha+1)},\quad  w_2=\frac{\Gamma(\alpha+1)}{c_2^{\alpha}\, \Gamma(2\alpha+1)},\quad
w_1=\frac{1}{\Gamma(\alpha+1)}-\frac{\Gamma(\alpha+1)}{c_2^{\alpha}\, \Gamma(2\alpha+1)}\,.
\end{equation}
inserting (\ref{EE14}) and (\ref{E18}) in (\ref{E13}) we get
\begin{align}
y_{n+1}&=y_n+\frac{h^{\alpha}}{\Gamma(\alpha+1)}f_n+\frac{h^{2\alpha}}{\Gamma(2\alpha+1)}\left( f_t^{\alpha}+f_nf_y\right)\nonumber \\
&\quad +\frac{c_2^{\alpha}h^{3\alpha}}{\Gamma(2\alpha+1)}
\left[ \frac{\Gamma(\alpha+1)}{\Gamma(2\alpha+1)}f_{t,t}^{\alpha,\alpha}+\frac{1}{2\Gamma(\alpha+1)}f_n^2f_{y,y}+
\frac{1}{\Gamma(\alpha+1)}f_nf_{t,y}^{\alpha,1}\right] +\cdots\label{E16}
\end{align}\label{Tn}
Subtracting (\ref{E16}) from (\ref{EE13}), we obtain the local
truncation error $T_n$

\begin{align}
T_n&=y(t_{n+1})-y_{n+1}=
h^{3\alpha}\left(\frac{1}{\Gamma(3\alpha+1)}-\frac{c_2^{\alpha}\Gamma(\alpha+1)}{(\Gamma(2\alpha+1))^2}   \right)  f_{t,t}^{\alpha,\alpha}\nonumber\\
&\hspace{2.8 cm}+ h^{3\alpha}   \left(
\frac{2}{\Gamma(3\alpha+1)}-\frac{c_2^{\alpha}}{2\Gamma(\alpha+1)}
\right)   f_n^2f_{y,y}
 \nonumber\\
 &\hspace{2.8 cm} +  h^{3\alpha}
 \left(  \frac{1}{\Gamma(3\alpha+1)}-\frac{c_2^{\alpha}}{\Gamma(\alpha+1)}  \right)  f_nf_{t,y}^{\alpha,1}+\cdots \label{tr1}
\end{align}

We conclude that no choice of the parameter $c_2^{\alpha}$ will make
the leading term of $T_n$ vanish for all functions $f(t,y)$.
Sometimes the free parameters are chosen to minimize the sum of the
absolute values of the coefficients in $T_n$. Such a choice is
called optimal choice. Obviously the minimum of $| T_n|$ occurs for
 $c_2^{\alpha}=\frac{(\Gamma(2\alpha+1))^2}{\Gamma(3\alpha+1)\Gamma(\alpha+1)}$,
 $c_2^{\alpha}=\frac{4\Gamma(\alpha+1)}{\Gamma(3\alpha+1)}$ or
  $c_2^{\alpha}=\frac{\Gamma(\alpha+1)}{\Gamma(3\alpha+1)}$.

From (\ref{tr1}) we have $\frac{T_n}{h^{\alpha}}=(h^{\alpha})^2$. So, we deduce that the 2-stage EFORK method (\ref{E10}) is of order ${2\alpha}$. \\

Now,  the 2-stage EFORK method by listing the coefficients is as
follows:
 \begin{table}[!h]
\begin{center}
\begin{tabular}{c | c c }
$c_2$     &\quad$a_{21}$\\
\hline &\quad $w_1$ &\quad $w_2$\\
 \end{tabular}\,,
\end{center}
\end{table}
%%%%%%%%%%%%%%%%%%%%%%%%%%%%
\begin{table}[!h]
\begin{center}
\begin{tabular}{c | c c }
$\left( \frac{2\Gamma(\alpha+1)^2}{\Gamma(2\alpha+1)}\right) ^\frac{1}{\alpha}$      &\quad$\frac{2\Gamma(\alpha+1)}{\Gamma(2\alpha+1)}$\\
\hline &\quad $\frac{1}{2\Gamma(\alpha +1)}$ &\quad $\frac{1}{2\Gamma(\alpha +1)}$\\
 \end{tabular}\,.
\end{center}
\end{table}
%%%%%%%%%%%%%%%%%%%%%%%%%%%
\\Also, the optimal cases of 2-stage EFORK method are:
\begin{table}[!h]
\begin{center}
%\caption{\small{ Optimal 2-stage method for $\alpha=1/2$, $T=1$ }}
\begin{tabular}{c | c c }
$\left( \frac{(\Gamma(2\alpha+1))^2}{\Gamma(3\alpha+1)\Gamma(\alpha+1)}\right) ^\frac{1}{\alpha}$      &\quad$\frac{(\Gamma(2\alpha+1))^2}{\Gamma(3\alpha+1)\Gamma(\alpha+1)^2}$\\
\hline &\quad        $\frac{1}{\Gamma(\alpha+1)}- \frac{\Gamma(3\alpha+1)\Gamma(\alpha+1)^2}{\Gamma(2\alpha+1)^3}$               &\quad  $\frac{\Gamma(3\alpha+1)\Gamma(\alpha+1)^2}{\Gamma(2\alpha+1)^3}$  \\
 \end{tabular}\,,
\end{center}
\end{table}
%%%%%%%%%%%%%%%%%%%%%%%%%%%%%%%%%
\begin{table}[!h]
\begin{center}
%\caption{\small{ Optimal 2-stage method for $\alpha=1/2$, $T=1$ }}
\begin{tabular}{c | c c }
$\left( \frac{4\Gamma(\alpha+1)}{\Gamma(3\alpha+1)}\right)^\frac{1}{\alpha} $     &\quad$\frac{4}{\Gamma(3\alpha+1)}$\\
\hline &\quad        $\frac{1}{\Gamma(\alpha+1)}- \frac{\Gamma(3\alpha+1)}{4\Gamma(2\alpha+1)}$  &\quad  $\frac{\Gamma(3\alpha+1)}{4\Gamma(2\alpha+1)}$  \\
 \end{tabular}\,,
\end{center}
\end{table}
%%%%%%%%%%%%%%%%%%%%%%%%%
\begin{table}[!h]
\begin{center}
\begin{tabular}{c | c c }
$\left( \frac{\Gamma(\alpha+1)}{\Gamma(3\alpha+1)}\right) ^\frac{1}{\alpha}$ &\quad$\frac{1}{\Gamma(3\alpha+1)}$\\
\hline &\quad $\frac{1}{\Gamma(\alpha+1)}- \frac{\Gamma(3\alpha+1)}{\Gamma(2\alpha+1)}$  &\quad  $\frac{\Gamma(3\alpha+1)}{\Gamma(2\alpha+1)}$  \\
 \end{tabular}\,.
\end{center}
\end{table}
%%%%%%%%%%%%%%%%%%%   3-stage
\subsection{ EFORK method of order $3 \alpha$}
Following (\ref{E8})-(\ref{E9}), we define a 3-stage EFORK method as
\begin{eqnarray}
&& K_1=h^{\alpha}f(t_n,y_n),\nonumber\\
&& K_2=h^{\alpha}f(t_+c_2h,y_n+a_{21}K_1),\nonumber\\
&& K_3=h^{\alpha}f(t_n+c_3h,y_n+a_{31}K_1+a_{32}K_2),\nonumber\\
&&y_{n+1}=y_n+w_1K_1+w_2K_2+w_3K_3\, .\label{E15}
\end{eqnarray}
where unknown parameters $\lbrace c_i \rbrace_{i=2}^3, \lbrace a_{ij}\rbrace_{i=2,j=1}^{3,i-1}$, and $\left\lbrace w_i\right\rbrace _{i=1}^3$
have to be determined accordingly. By using the same procedure as we
did for 2-stage EFORK method, expanding $K_1$, $K_2$ and $K_3$,
comparing with (\ref{EE13}) and matching coefficients of powers of
$h^{\alpha}$, we obtain the following equations:
\begin{eqnarray}\label{E17}
&& w_1+w_2+w_3=\frac{1}{\Gamma(\alpha+1)},  \qquad a_{21}=\frac{c_2^{\alpha}}{\Gamma(\alpha+1)},\nonumber\\
&&w_2c_2^{2\alpha}+w_3c_3^{2\alpha}=\frac{\Gamma(2\alpha+1)}{\Gamma(3\alpha+1)}, \qquad a_{31}+a_{32}=\frac{c_3^{\alpha}}{\Gamma(\alpha+1)},   \nonumber\\
&&w_2c_2^{\alpha}+w_3c_3^{\alpha}=\frac{\Gamma(\alpha+1)}{\Gamma(2\alpha+1)}, \qquad w_3a_{32}c_2^{\alpha}=\frac{\Gamma(\alpha+1)}{\Gamma(3\alpha+1)},\nonumber\\
\end{eqnarray}
Now, we have six equations with eight unknown parameters. According
to Butcher tableau for 3-stage EFORK method, we have:
 \begin{table}[!h]
\begin{center}
\begin{tabular}{c |c c c}
$c_2$     &\quad$a_{21}$\\
$c_3$      &\quad$a_{31}$   &\quad$a_{32}$\\
\hline &\quad $w_1$ &\quad $w_2$ &\quad $w_3$\\
 \end{tabular}
\end{center}
\end{table}
%%%%%%%%%%%%%%%%%%%%%%%%%%%%%%%%
\\ If $c_2$ and $c_3$ are arbitrarily chosen, we calculate weights $\left\lbrace w_i\right\rbrace _{i=1}^3$ and coefficients
$\lbrace a_{ij}\rbrace_{i=2,j=1}^{3,i-1}$ from (\ref{E17}) as:

\begin{table}[!h]
\begin{center}
\begin{tabular}{c |c c c}
$\left( \frac{1}{2\alpha !}\right) ^\frac{1}{\alpha}$ & $\frac{1}{2(\alpha !)^2}$\\
$\left( \frac{1}{4\alpha !}\right) ^\frac{1}{\alpha}$ &
$\frac{(\alpha !)^2(2\alpha)!+
2(2\alpha)!^2-(3\alpha)!}{4(\alpha !)^2(2(2\alpha)!^2-(3\alpha)!)}$   &  $-\frac{(2\alpha)!}{4(2(2\alpha)!)-(3\alpha)!}$\\
\hline &   $ \frac{8(\alpha
!)^2(2\alpha)!}{(3\alpha)!}-\frac{6(\alpha
!)^2}{(2\alpha)!}+\frac{1}{\alpha !} $ &   $ \frac{2(\alpha
!)^2(4(2\alpha)!^2-(3\alpha)! )}{(2\alpha)!(3\alpha)!} $
 &   $- \frac{8(\alpha !)^2(2(2\alpha)!^2-(3\alpha)! )}{(2\alpha)!(3\alpha)!} $\\
 \end{tabular}
\end{center}
\end{table}

As a result, we obtain $\frac{T_n}{h^{\alpha}}=(h^{\alpha})^3$. In a
similar procedure with 2 and 3 stages EFORK methods we can construct
s-stages EFORK methods for $s>3.$

% {One can see that the following
%relation can be obtain by using above procedure for
%$s\geq4$.}\\\\
%
%

As we can see, to obtain the higher fractional order Runge-Kutta methods, we
must consider a method with additional stages. In next section, we
express implicit  fractional order Runge-Kutta (IFORK) methods with low stages and
high orders.

%%%%%%%%%%%%%%%%%%%%%%%%%%%%
%%%%%%%%%%%%%%%%%%%%%%%%
%%%%%%%%%%%%%%%%%%%%%%%

\section{IFORK methods}\label{se2}
We define a s-stage IFORK method by the following equations:
\begin{equation}\label{ims}
K_i=\frac{1}{s}h^{\alpha}\sum_{k=1}^s f(t_n+c_{ik}h,y_n+\sum_
{j=1}^sa_{ij}K_j)
 ,\quad i=1,2,...,s,
\end{equation}
and
\begin{equation}\label{11}
y_{n+1}=y_n+\sum_{i=1}^sw_iK_i,
\end{equation}
where
\begin{equation}\label{22}
\frac{c_{i1}^{\alpha}+c_{i2}^{\alpha}+...+c_{is}^{\alpha}}{\alpha
!}=s(a_{i1}+a_{i2}+...+a_{is}), \quad i=1,2,...,s
\end{equation}
and the parameters $\left\lbrace a_{ij}\right\rbrace_{i,j=1}^{s,s} $,
$\left\lbrace w_i \right\rbrace_{i=1}^s $ are arbitrary. We state
the IFORK method by listing the coefficients as follows:
%%%%%%%%%%%%%%%%%%%%%%%%%%
 \begin{table}[!h]
\begin{center}
\begin{tabular}{c c c c|c c c c}
$c_{11}$&\quad$c_{12}$ &\quad $\cdots$&\quad$c_{1s}$&\quad$a_{11}$ &\quad$a_{12}$&\quad $\cdots$&\quad$a_{1s}$\\
$c_{21}$ &\quad $c_{22}$ &\quad $\cdots$&\quad$c_{2s}$&\quad$a_{21}$ &\quad$a_{22}$&\quad $\cdots$&\quad$a_{2s}$\\
$\vdots$&\quad$\vdots$&\quad$\vdots$&\quad$\vdots$&\quad$\vdots$&\quad$\vdots$&\quad$\vdots$&\quad$\vdots$\\
$c_{s1}$&\quad$c_{s2}$&\quad $\cdots$ &\quad$c_{ss}$&\quad$a_{s1}$ &\quad$a_{s2}$&\quad $\cdots$&\quad$a_{ss}$\\
\hline &\quad &\quad &\quad  &\quad$w_1$ &\quad $w_2$&\quad$\cdots$ &\quad$w_s$\\
\end{tabular}
\end{center}
\end{table}
%%%%%%%%%%%%%%%%%%%%%%%%%%%%%%%%%%%%%%%%
\\Since the  functions $K_i$ are  defined  by  a  set of $s$ implicit
equations, the  derivation  of  the  implicit methods is
complicated. Therefore, without loss of generality, only the case
$s=2$ is investigated.

Consider (\ref{ims})-(\ref{22}) with $s=2$ as
\begin{equation}\label{iy1}
K_i=\frac{1}{2} h^{\alpha}\left[ f(t_n+c_{i1}h,y_n+a_{i1}K_1+a_{i2}K_2)+
f(t_n+c_{i2}h,y_n+a_{i1}K_1+a_{i2}K_2)\right] ,\quad i=1,2
\end{equation}
\begin{equation}\label{y1}
y_{n+1}=y_n+w_1K_1+w_2K_2
\end{equation}
where
\begin{equation}\label{iy2}
\frac{c_{i1}^{\alpha}+c_{i2}^{\alpha}}{\alpha !}=2(a_{i1}+a_{i2}),
\quad i=1,2
\end{equation}
By using the similar procedure as we did for EFORK method, we
expand $K_i$ about the point $(t_n,y_n)$, where we apply Caputo
Taylor formula  (\ref{EE13}) about $t_n$ and standard integer order
Taylor formula about $y_n$.
\begin{align}
K_i&=\frac{1}{2} h^{\alpha}[ 2
f_n+\frac{(c_{i1}^{\alpha}+c_{i2}^{\alpha})h^{\alpha}}{\alpha
!}f_t^{\alpha}+2(a_{i1}K_1+a_{i2}K_2)f_y+
\frac{(c_{i1}^{\alpha}+c_{i2}^{\alpha})h^{2\alpha}}{(2\alpha) !}f_{t,t}^{\alpha,\alpha}\nonumber\\
&\quad  +{(a_{i1}K_1+a_{i2}K_2)^2}f_{y,y}+\frac{(c_{i1}^{\alpha}+c_{i2}^{\alpha})h^{\alpha}}{\alpha
!}(a_{i1}K_1+a_{i2}K_2)f_{t,y}^{\alpha,1}
\nonumber\\
&\quad   +\frac{(c_{i1}^{3\alpha}+c_{i2}^{3\alpha})h^{3\alpha}}{(3\alpha) !}f_{t,t,t}^{\alpha,\alpha,\alpha}+
\frac{(c_{i1}^{2\alpha}+c_{i2}^{2\alpha})h^{2\alpha}}{(2\alpha)
!}(a_{i1}K_1+a_{i2}K_2)f_{t,t,y}^{\alpha,\alpha,1}
 \nonumber\\
&\quad   +\frac{(c_{i1}^{\alpha}+c_{i2}^{\alpha})h^{\alpha}}{\alpha !}\frac{(a_{i1}K_1+a_{i2}K_2)^2}{2}f_{t,y,y}^{\alpha,1,1}+\frac{(a_{i1}K_1+a_{i2}K_2)^3}{3}f_{y,y,y}+\cdots]   , \label{im1}
\end{align}
where $i=1,2$\,.

Since equations (\ref{im1}) are implicit, we cannot obtain the
explicit forms for $K_1$ and $K_2$. To determine the explicit form
$K_i$, we consider
\begin{equation}\label{E19}
K_i=h^{\alpha}A_i+h^{2\alpha}B_i+h^{3\alpha}C_i+\cdots,\quad i=1,2
\end{equation}
where $A_i$, $B_i$ and $C_i$ are unknowns. Substituting (\ref{E19})
into (\ref{im1}) and matching the coefficients of powers of
$h^{\alpha}$, we get
\begin{align}
A_i&=f_n,\nonumber\\
B_i&=\frac{c_{i1}^{\alpha}+c_{i2}^{\alpha}}{2(\alpha !)}f_t^{\alpha}+a_{i1}ff_y+a_{i2}ff_y=\frac{c_{i1}^{\alpha}+c_{i2}^{\alpha}}{2(\alpha !)}D^{\alpha}f,\nonumber\\
C_i&=\left( a_{i1}\frac{c_{11}^{\alpha}+c_{12}^{\alpha}}{2(\alpha
!)}+a_{i2}\frac{c_{21}^{\alpha}+c_{22}^{\alpha}}{2(\alpha !)}  \right)
f_y D^{\alpha} f+
\frac{c_{i1}^{2\alpha}+c_{i2}^{2\alpha}}{2(2\alpha)!}f_{t,t}^{\alpha,\alpha}
\nonumber\\
  &\quad+\frac{1}{4}\left( \frac{c_{i1}^{\alpha}+c_{i2}^{\alpha}}{\alpha
!}\right) ^2\left( \frac{1}{2} f^2f_{yy}
+ff_{t,y}^{\alpha,1}\right), \nonumber\\
\vdots
 \label{E20}
\end{align}
Inserting (\ref{E19}) and (\ref{E20}) into (\ref{y1}), we have
 \begin{align}
 y_{n+1}&=y_n+h^{\alpha}\left[ w_1A_1+w_2A_2\right]+h^{2\alpha}\left[  w_1B_1+w_2B_2  \right] +h^{3\alpha}\left[  w_1C_1+w_2C_2  \right] +\cdots\label{y2}
  \end{align}
  Comparing (\ref{y2}) with (\ref{EE13}) and equating the coefficient of powers of $h^{\alpha}$, we can get IFORK method
of different orders.

\subsection{IFORK method of order $2 \alpha$}
To obtain an IFORK method of order $2\alpha$, we equate the
coefficients of $h^{\alpha}$ and $h^{2\alpha}$ in (\ref{EE13}) and
(\ref{y2}) correspondingly, to get
 \begin{eqnarray*}
 &&w_1+w_2=\frac{1}{\alpha !},\\
 &&w_1\frac{c_{11}^{\alpha}+c_{12}^{\alpha}}{\alpha !}+w_2\frac{c_{21}^{\alpha}+c_{22}^{\alpha}}{\alpha !}=\frac{2}{(2\alpha )!},
 \end{eqnarray*}
 where
 \begin{equation*}
 2(a_{11}+a_{12})=\frac{c_{11}^{\alpha}+c_{12}^{\alpha}}{\alpha !}, \quad    2(a_{21}+a_{22})=\frac{c_{21}^{\alpha}+c_{22}^{\alpha}}{\alpha !}.
 \end{equation*}
 There are now six arbitrary parameters to be prescribed. If we neglect $K_2$, i.e, if we choose $a_{21}=a_{22}=a_{12}=0$, $w_2=0$, from the above equations, we find
\begin{equation*}
w_1=\frac{1}{\alpha !}, \quad
c_{11}^{\alpha}+c_{12}^{\alpha}=\frac{2(\alpha !)^2}{(2\alpha)!},
\quad a_{11}=\frac{\alpha !}{(2\alpha)!}.
\end{equation*}
%%%%%%%%%%%%%%%%%%%%%%%%%%%%%%%%%%%%%%%
Therefor, a 1-stage IFORK method of order $2\alpha$ is obtained as
follows:
 \begin{eqnarray}
&& K_1=\frac{1}{2} h^{\alpha}\left[ f(t_n+c_{11}h,y_n+a_{11}K_1)+f(t_n+c_{12}h,y_n+a_{11}K_1)\right] ,\nonumber\\
&&y_{n+1}=y_n+w_1K_1\label{1i}.
 \end{eqnarray}

\subsection{IFORK method of order $3 \alpha$}
 Also, we can get IFORK method of order $3\alpha$ with 2-stage (\ref{iy1})-(\ref{iy2}), when
 equating the coefficients of $h^{\alpha}$, $h^{2\alpha}$ and $h^{3\alpha}$ in (\ref{EE13}) and
(\ref{y2}) accordingly. In such case, we obtain the following system
of equations
 \begin{eqnarray*}
 &&w_1+w_2=\frac{1}{\alpha !},\\
 &&w_1\frac{c_{11}^{\alpha}+c_{12}^{\alpha}}{\alpha !}+w_2\frac{c_{21}^{\alpha}+c_{22}^{\alpha}}{\alpha !}=\frac{2}{(2\alpha )!},\\
 &&w_1 \left(  a_{11}\frac{c_{11}^{\alpha}+c_{12}^{\alpha}}{\alpha !}+a_{12}\frac{c_{21}^{\alpha}+c_{22}^{\alpha}}{\alpha !}
  \right) +w_2\left(  a_{21}\frac{c_{11}^{\alpha}+c_{12}^{\alpha}}{\alpha !}+a_{22}\frac{c_{21}^{\alpha}+c_{22}^{\alpha}}{\alpha !}                 \right) =\frac{2}{(3\alpha)!},\\
 && w_1 \frac{c_{11}^{2\alpha}+c_{12}^{2\alpha}}{(2\alpha) !}+w_2\frac{c_{21}^{2\alpha}+c_{22}^{2\alpha}}{(2\alpha) !}                =\frac{2}{(3\alpha)!},\\
  && w_1 \left( \frac{c_{11}^{\alpha}+c_{12}^{\alpha}}{\alpha !}\right) ^2+w_2\left( \frac{c_{21}^{\alpha}+c_{22}^{\alpha}}{\alpha !}\right) ^2=\frac{8}{(3\alpha)!},
 \end{eqnarray*}
 where
 \begin{equation*}
 2(a_{11}+a_{12})=\frac{c_{11}^{\alpha}+c_{12}^{\alpha}}{\alpha !}, \quad    2(a_{21}+a_{22})=\frac{c_{21}^{\alpha}+c_{22}^{\alpha}}{\alpha !}.
 \end{equation*}
 The three free parameters can be chosen in such a way that $K_1$ or $K_2$ be explicit. If we want $K_1$ to be explicit, we choose
 $$c_{11}=a_{11}=a_{12}=0 .   $$

 Thus, the IFORK method of order $3\alpha$ which is explicit in $K_1$ is given by
 \begin{eqnarray}
 &&K_1=\frac{1}{2} h^{\alpha}\left[ f(t_n,y_n)+ f(t_n+c_{12}h,y_n)\right] ,\nonumber\\
 && K_2=\frac{1}{2} h^{\alpha}\left[ f(t_n+c_{21}h,y_n+a_{21}K_1+a_{22}K_2)+f(t_n+c_{22}h,y_n+a_{21}K_1+a_{22}K_2)\right] ,\nonumber\\
 && y_{n+1}=y_n+w_1K_1+w_2K_2\, .\label{2i}
 \end{eqnarray}
%%%%%%%%%%%%%%%%%%%%%%%%%%%%%%%%%%
Now, we can write the IFORK method of order $3\alpha$ as:
 %%%%%%%%%%%%%%%%%%%%%%%%%%%%%%
\begin{table}[!h]
\begin{center}
\begin{tabular}{c c|c c }
$c_{11}$&\quad$c_{12}$ &\quad$a_{11}$ &\quad$a_{12}$\\
$c_{21}$ &\quad $c_{22}$ &\quad$a_{21}$ &\quad$a_{22}$\\
\hline &\quad &\quad $w_1$ &\quad $w_2$
\end{tabular},
\end{center}
\end{table}
%%%%%%%%%%%%%%%%%%%%%%%

\begin{table}[!h]
\begin{center}
\begin{tabular}{c c|c c }
$0$&$0$ &$0$ &$0$\\
$\left(\frac{(2\alpha)!(2(\alpha!) -\sqrt{2}\sqrt{ (2\alpha)!
-2(\alpha!)^2 })}{(3\alpha)!}\right) ^{\frac{1}{\alpha}}$ &
$\left(\frac{(2\alpha)!(2(\alpha!) + \sqrt{2}\sqrt{ (2\alpha)!
-2(\alpha!)^2  })}{(3\alpha)!}\right)
^{\frac{1}{\alpha}}$    &     $\frac{(2\alpha)!}{(3\alpha)!}$ &\hspace{-.5 cm}  $\frac{(2\alpha)!}{(3\alpha)!}$\\
\hline & & $  \frac{1}{\alpha!}-\frac{(3\alpha)!}{2(2\alpha)!^2}  $   &  $\frac{(3\alpha)!}{2(2\alpha)!^2}$ \\
\end{tabular}
\end{center}
\end{table}
%%%%%%%%%%%%%%%%%%%%%%%%%%%%%%%%%%%%%%%%%%%%%%%%%%%%%%%%%%
%%%%%%%%%%%%%%%%%%%%% Thoeretical analysis%%%%%%%%%%%%%%%%%%%%%%%%%%%%%%%%%%%%%%%%%%%%%%%%%%%%%%%%%%%%%%%%
\section{Theoretical analysis}\label{st}
The obtained difference approximation to the FDEs in FORK methods, does not guarantee  that  the
solution of the difference equation  can approximates  the exact
solution of the FDE correctly. Here the
convergence analysis of the FORK methods which arises from
some conditions for which the difference solutions converge  to the
exact solution is investigated.

In this section, firstly we consider a definition of consistency of
the discussed methods  in section \ref{se1} and \ref{se2}.
%\section{ Theoretical analysis}\label{st}
%In this section after expressing the existence and uniqueness of
%solution of the fractional differential equation (\ref{E6}) from
%\cite{17}, we also discuss the consistency and the stability, as two
%major milestones of convergence, of RK methods and consider a new
%definition of consistency which is used to the methods discussed in
%section \ref{se1}.
%\begin{thm}
%\cite{17}.Let $\alpha>0$, $y_0\in \mathbb{R}$, $K>0$ and $h^*>0$.
%and let the function $f:G\rightarrow \mathbb{R}$ be continuous and
%fulfil a Lipschitz condition with respect to the second variable,
%i.e.
%$$ |f(t,y_1)-f(t,y_2)|\leq L|y_1-y_2|$$
%with some constant $L>0$ independent of $t, y_1$, and $y_2$. Define
%$G=\left\lbrace (t,y):t\in [0,h^*], |y-y_0|\leq K  \right\rbrace $,
%$M=Sup_{(t,z)\in G}|f(t,z)|$ and
%\[h = \left\{ \begin{gathered}
%  \begin{array}{*{20}{c}}
%  {h^*}&{}&{if}&{M = 0,}
%\end{array} \hfill \\
%  \begin{array}{*{20}{c}}
%  {\min \{ h^*,{{(K\Gamma (\alpha + 1)/M)}^{1/\alpha}}\} }&{}&{else.}&{}
%\end{array} \hfill \\
%\end{gathered}  \right.\]
%Then, there exists a uniquely function  $y\in C[0,h]$solving the
%initial-value problem (\ref{E6}).
%\end{thm}
%%%%%%%%%%%%%%%%%%%%%%%%%%%%%%%%%%%%%%%%%%%%%%%%%%%%%%%
\subsection{ Consistency}
The EFORK and IFORK methods considered before belong to the class of
methods which are characterized by the use of $y_n$ on the
computation of $y_{n+1}$. These family of one-step methods admits
the following representation
%In order to investigate the consistency of the FORK method, we
%introduce a new definition of consistency as follows:
%\begin{defn}\label{D1}
%The one-step method (\ref{C1}) is said to be consistent if
%\begin{equation*}
%\Phi(t,y,0)=\frac{1}{\Gamma(\alpha+1)}f(t,y).
%\end{equation*}
%where $\Phi(t_n,y_n,h)$ is a functions of arguments $t_n$, $y_n$,
%$h$ and in addition depends on the right hand side of (\ref{E6}).
%\end{defn}

\begin{eqnarray}\label{cons1}
&&y_{n+1}=y_n+h^{\alpha} \Phi(t_n,y_n,y_{n+1},h),~~~~
n=0,...,N^m-1,\\
&& y_0 = y(t_0).\nonumber
\end{eqnarray}
where $\Phi : [t_0, T]\times\mathbb{R}^2\times[0,h_0]\to\mathbb{R}$
and for the particular case of the explicit methods we have the
representation
\begin{eqnarray}\label{cons2}
&&y_{n+1}=y_n+h^{\alpha} \Phi(t_n,y_n,h),~~~~
n=0,...,N^m-1,\\
&& y_0 = y(t_0).\nonumber
\end{eqnarray}
with  $\Phi :
[t_0,T]\times\mathbb{R}\times[0,h_0]\rightarrow\mathbb{R}$.\\ We
define the truncation error $\tau_n$ by

\begin{equation}\label{der}
\tau_n=\frac{y_{n+1}-y_n}{h^\alpha}- \Phi(t_n,y_n,y_{n+1},h),~~~~\\
\end{equation}
The one-step method (\ref{cons1}) and (\ref{cons2}) is said
consistent with the equation (\ref{E6}), if
$\lim_{n\rightarrow\infty} \tau_n= 0, ~~ h N^m =T-t_0.$

Therefore from (\ref{der}) and (\ref{EE13}) we may write
\begin{eqnarray*}
&&\displaystyle{\lim_{h\to 0}}\tau_n=\displaystyle{\lim_{h\to
0}}\frac{y_{n+1}-y_n}{h^\alpha}- \displaystyle{\lim_{h\to
0}}\Phi(t_n,y_n,y_{n+1},h),\\
&&~~~~~~~~=\frac{1}{\Gamma(\alpha+1)}{}_{t_n}^cD_t^{\alpha}y(t_n)-
\displaystyle{\lim_{h\to 0}}\Phi(t_n,y(t_n),y(t_{n+1}),h),\\
&&~~~~~~~~=\frac{1}{\Gamma(\alpha+1)}F_n(t_n,y(t_n))-
\displaystyle{\lim_{h\to 0}}\Phi(t_n,y(t_n),y(t_n+h),h).
\end{eqnarray*}

Now we conclude that the proposed one-step EFORK and IFORK methods
are consistent if and only if
\begin{equation*}
\Phi(t,y,y,0)=\frac{1}{\Gamma(\alpha+1)}F_n(t,y),
\end{equation*}
or briefly
\begin{equation*}
\Phi(t,y,y,0)=\frac{1}{\Gamma(\alpha+1)}f(t,y).
\end{equation*}
Also in a similar manner, for explicit methods we have
\begin{equation*}
\Phi(t,y,0)=\frac{1}{\Gamma(\alpha+1)}f(t,y).
\end{equation*}

As an example, consider the 2-stage EFORK method (\ref{E10}) with
(\ref{EE14}),
\begin{equation*}
y_{n+1}=y_n+h^{\alpha}\left[
w_1f(t_n,y_n)+w_2f(t_n+c_2h,y_n+a_{21}K_1)  \right],
\end{equation*}
where in comparison to (\ref{cons2}), we have
\begin{equation*}
\Phi(t,y,h)=\left[ w_1f(t,y)+w_2f(t+c_2h,y+a_{21}K_1)  \right],
\end{equation*}
or when h tends to $0 $ yields
\begin{equation*}
\Phi(t,y,0)=(w_1+w_2)f(t,y),
\end{equation*}
and by using (\ref{EE14}), we may write
$$\Phi(t,y,0)=\frac{1}{\Gamma(\alpha+1)}f(t,y). $$
Therefor, 2-stage EFORK method (\ref{E10}) is consistent. Also,
3-stage EFORK method (\ref{E15}) is consistent, according to
\begin{eqnarray*}
&&\Phi(t,y,h)=\left[ w_1f(t,y)+w_2f(t+c_2h,y+a_{21}K_1) +w_3f(t+c_3h,y+a_{31}K_1+a_{32}K_2) \right],\\
&&\Phi(t,y,0)=(w_1+w_2+w_3)f(t,y),\\
\end{eqnarray*}
and so from (\ref{E17})
$$\Phi(t,y,0)=\frac{1}{\Gamma(\alpha+1)}f(t,y). $$

Similarly, we can show the consistency of all proposed FORK methods
in sections \ref{se1} and \ref{se2}.
%%%%%%%%%%%%%%%%%%%%%%%%%%%%%%%%%%%%%%%
\subsection{Convergence analysis}
In this section we investigate the convergence behavior of the
proposed FORK methods (without loss of generality we consider only
explicit FORK  methods). To do so, we express a definition of
regularity from \cite{4}.
\begin{defn}
A one-step method of the form (\ref{cons2})
\begin{equation}\label{C1}
y_{n+1}=y_n+h^{\alpha}\Phi(t_n,y_n,h),\quad n=0,1,2,...,N^m-1,
\end{equation}
is said to be regular if the function $\Phi(t,y,h)$ is defined and continuous in domain $t\in[0,T]$, $y\in[0,T^*]$ and $h\in[0,h_0]$ ($h_0$ is a positive constant) and if there exists a constant $L$ such that
\begin{equation*}
|\Phi(t,y,h)-\Phi(t,z,h)|\leq L|y-z|,
\end{equation*}
for every $t\in[0,T]$, $y,z\in[0,T^*]$ and $h\in[0,h_0]$.
\end{defn}
%%%%%%%%%%%%%%%%%%%%%%%%%%%
To discuss the convergence of the EFORK methods, at first we prove
the given methods in section $3$ are regular. We know from Theorem
\ref{Tu}, that $f(t,y)$ satisfies a Lipschitz condition with respect
to second variable. Thus
\begin{align}
 |\Phi(t_n,y_n,h)-\Phi(t_n,y_n^*,h)|&=h^{-\alpha}\vert\sum_{i=1}^{s}w_iK_i-\sum_{i=1}^{s}w_iK_i^*\vert \nonumber\\
 &\leq h^{-\alpha}\left( w_1\vert k_1-k_1^*\vert +w_2 \vert k_2-k_2^*\vert  +\ldots +w_s\vert k_s-k_s^*\vert \right) \nonumber\\
 & \leq w_1L|y_n-y_n^*|+w_2L|y_n+a_{21}k_1-y_n^*-a_{21}k_1^*|+\ldots \nonumber\\
 &+w_sL|y_n+a_{s1}k_1+a_{s2}k_2+\ldots+a_{ss}k_s-y_n^*-a_{s1}k_1^*\nonumber\\
 &-a_{s2}k_2^*-\ldots-a_{ss}k_s^*| \nonumber\\
 &\leq \ldots \leq L^*|y_n-y_n^*|.\nonumber
\end{align}
Therefor the function $\Phi$ satisfies a Lipschitz condition in $y$
and it is also continuous, thus EFORK methods are regular. To
establish the convergence behavior, we need following  Lemma from
\cite{4}.
\begin{lem}\label{lem}
Let $\omega_0, \omega_1,\omega_2,...$ be a sequence of real positive
numbers which satisfy
$$\omega_{n+1}\leq (1+\zeta)\omega_n+\mu, \qquad n=0,1,2... $$
where $\zeta$, $\mu$ are positive constants. Then
\begin{equation*}
 \omega_n\leq e^{n\zeta}\omega_0+\left( \frac{e^{n\zeta}-1}{\zeta}\right)\mu, \qquad n=0,1,2,... \,.
\end{equation*}
\end{lem}
We now discuss on the behavior of the error $e_n=y(t_n)-y_n$ in
EFORK method for the initial-value problem (\ref{E6}).
%%%%%%%%%%%%%%%%%%%%%%%%%%%%%%%%%%%%%%%%%%%%%%%%%
\begin{thm}Consider the initial value problem (\ref{E6}) and
Let $f(t,y(t))$ is continuous and satisfy a Lipschitz condition with
Lipschitz constant $L$ and also
 $({}_{t_0}^cD_t^{\alpha})^{(s+1)}y(t)$ is continuous for $t\in [t_0,T]$,
Then the given EFORK method in section $3$ is convergent for
$m\alpha \geq 1$, if and only if it is consistent.
\end{thm}
\begin{proof}
Let EFORK method be consistent and the method can be written in the
form
\begin{equation}\label{sc1}
y_{n+1}=y_n+h^{\alpha}\Phi(t_n,y_n,h).
\end{equation}
The exact value $y(t_n)$ will satisfy
\begin{equation}\label{sc2}
y(t_{n+1})=y(t_n)+h^{\alpha}\Phi(t_n,y(t_n),h)+T_n,
\end{equation}
where $T_n$ is the truncation error. By subtracting (\ref{sc1}) from
(\ref{sc2}), we have
\begin{align*}
&|e_{n+1}| \leq | e_n|+h^{\alpha}|\left(
\Phi(t_n,y(t_n),h)-\Phi(t_n,y_n,h) \right)|+|T_n|.
\end{align*}
Now from regularity of the EFORK method it follows that
\begin{align*}
&|e_{n+1}|\leq |e_n|+h^{\alpha}L|y(t_n)-y_n|+|T_n|\\
& \;\qquad \;\leq (1+h^{\alpha}L)|e_n|+|T_n|.
\end{align*}
 By using the Lemma \ref{lem}, we have
 \begin{equation*}
 |e_{n} |\leq  (1+h^{\alpha}L)^n|e_0|+\left( \frac{e^{nh^{\alpha}L}-1}{h^{\alpha}L}\right) |T_n|,
 \end{equation*}
 where we assumed that the local truncation error for sufficiently large n is constant, i.e. $\mathcal{T}=T_n, n=0,1,2,...$. Also, assume that
 $e_0=0$ and $ |T_n|=O(h^{p \alpha}), p\geq 3$, therefor
 \begin{align*}
& |e_{n} |\leq O(h^{p\alpha})\left( \frac{e^{nh^{\alpha}L}-1}{h^{\alpha}L}\right).\\
\end{align*}

In section $3$ we assumed $N^m=\frac{T-t_0}{h},$ so  we have
 \begin{align*}
& |e_{n} |\leq O(h^{(p-1)\alpha})\left(
\frac{e^{(T-t_0)^{\frac{1}{m}}Lh^{\alpha-\frac{1}{m}}}-1}{L}\right).
\end{align*}

Thus, the EFORK methods of subsections $3.1-3.2$ are
convergent if $\alpha-\frac{1}{m}\geq 0$, i.e $m\alpha\geq 1$.\\

%%%%%%%%%%%%%%%%%%%%%
Now conversely, let EFORK method be convergent. It is sufficient we
give a limit of (\ref{sc1}) as $h$ tends to $0$. Now the proof of
theorem is complete.\\
\end{proof}
%%%%%%%%%%%%%%%%%%%%%%%%%%%%%%%%%%
\subsection{Stability analysis}
For stability analysis of the proposed method in section
\ref{se1}-\ref{se2}, we consider the FDE
\begin{eqnarray}
&& {}_{t_0}^cD_t^{\alpha} y(t)=\lambda y(t),\quad \lambda \in \mathbb{C},\quad 0<\alpha\leq1 ,        \nonumber\\
&& y(t_0)=y_0\label{s1}\,.
\end{eqnarray}
According to \cite{1}, the exact solution of (\ref{s1}) is
$y(t)=E_{\alpha}(\lambda (t-t_0)^{\alpha})y_0$. When
$\Re(\lambda)<0$, the solution of (\ref{s1}) asymptotically tends to
$0$ as $t\rightarrow \infty$.

\begin{figure}
\begin{center}
\includegraphics[width=5.5cm]{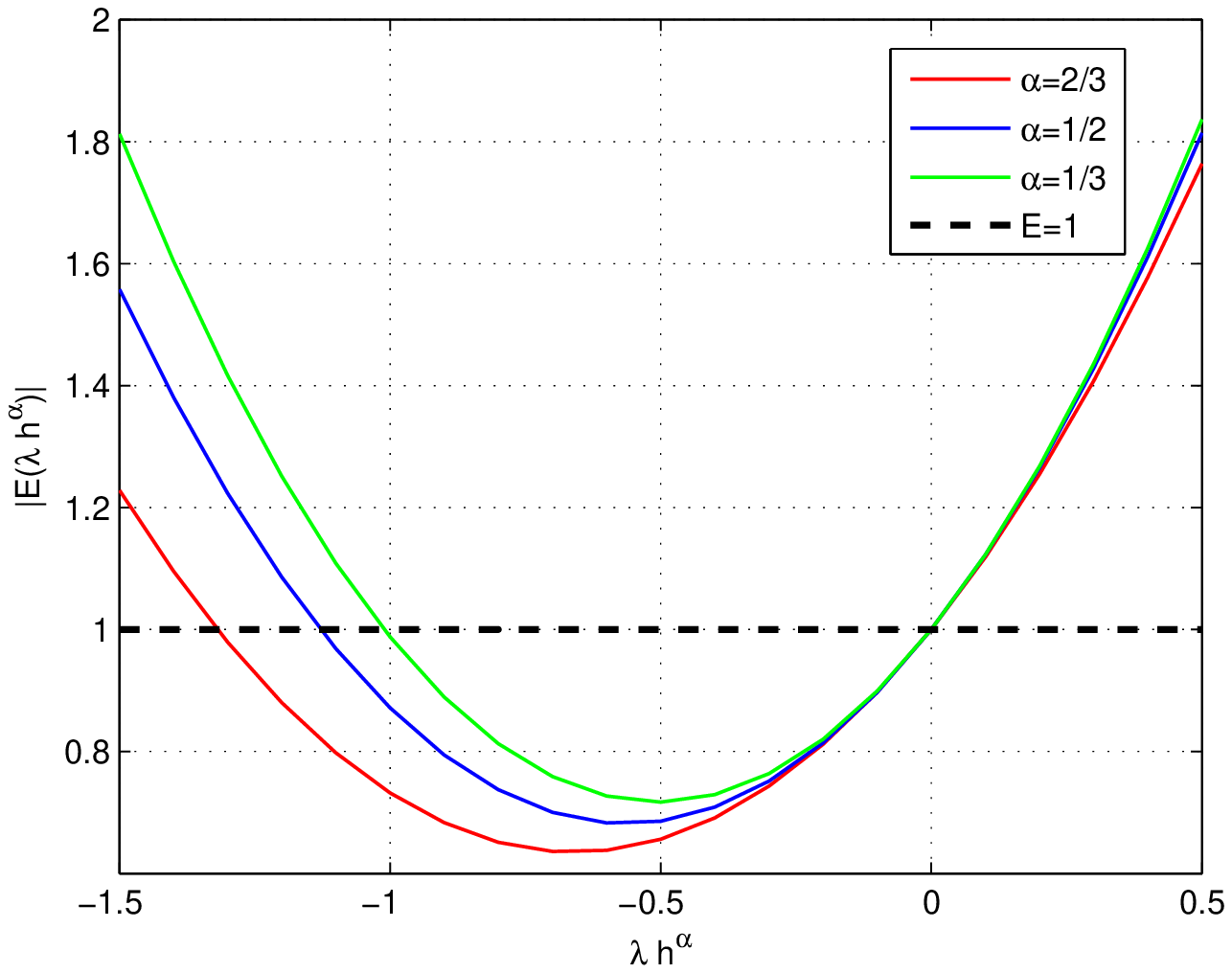}
\includegraphics[width=5.5cm]{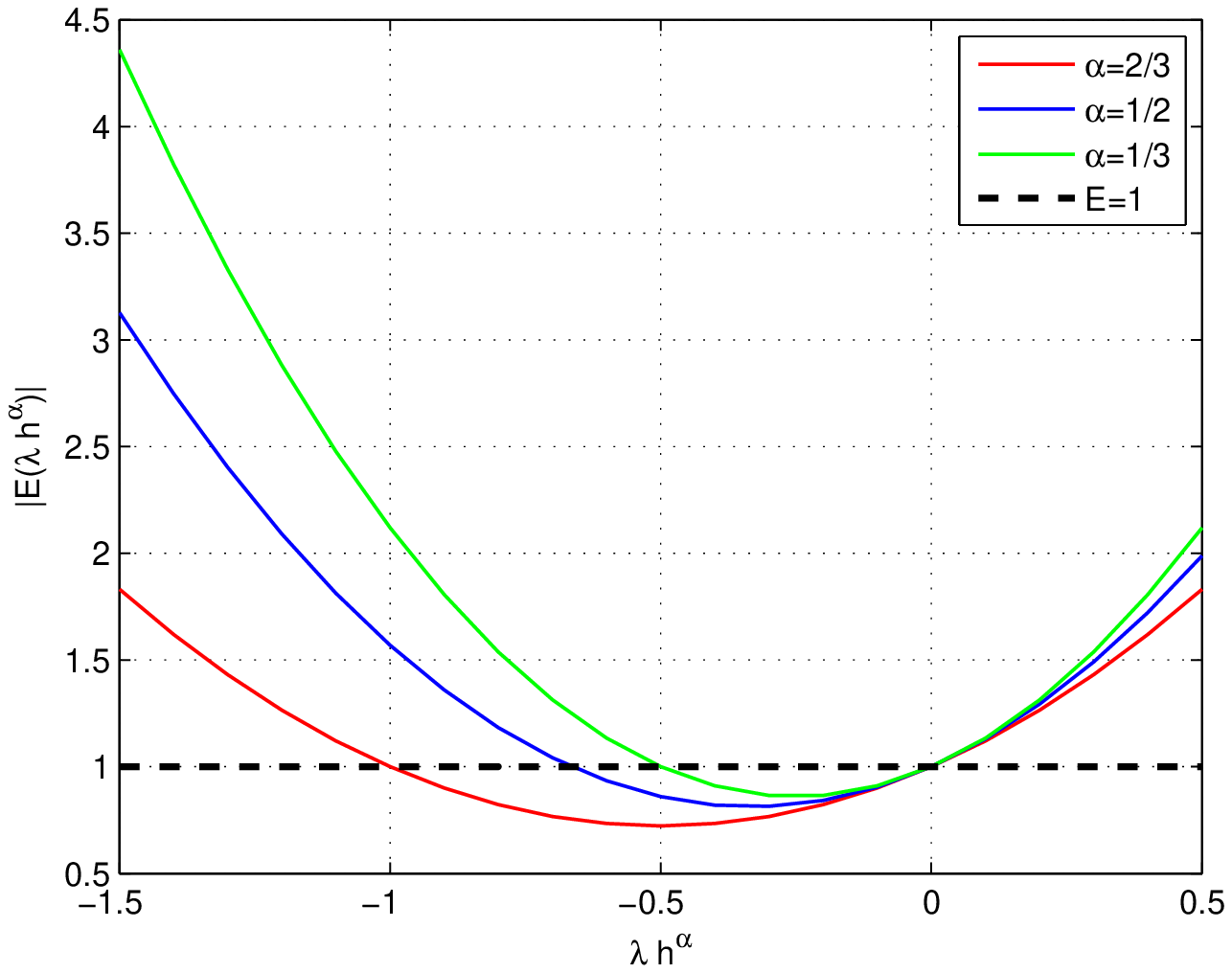}
\caption{\small{The graph of $E(\lambda h^{\alpha})$ for 2-stage
EFORK method (\ref{E10}) with $c_2^{\alpha}=\frac{2(\alpha
!)^2}{(2\alpha)!}$ ({\it left}), and
$c_2^{\alpha}=\frac{(\Gamma(2\alpha+1))^2}{\Gamma(3\alpha+1)\Gamma(\alpha+1)}$
({\it right}). }}
 \label{ps-1}
 \end{center}
\end{figure}

\begin{figure}
\begin{center}
\includegraphics[width=5.5cm]{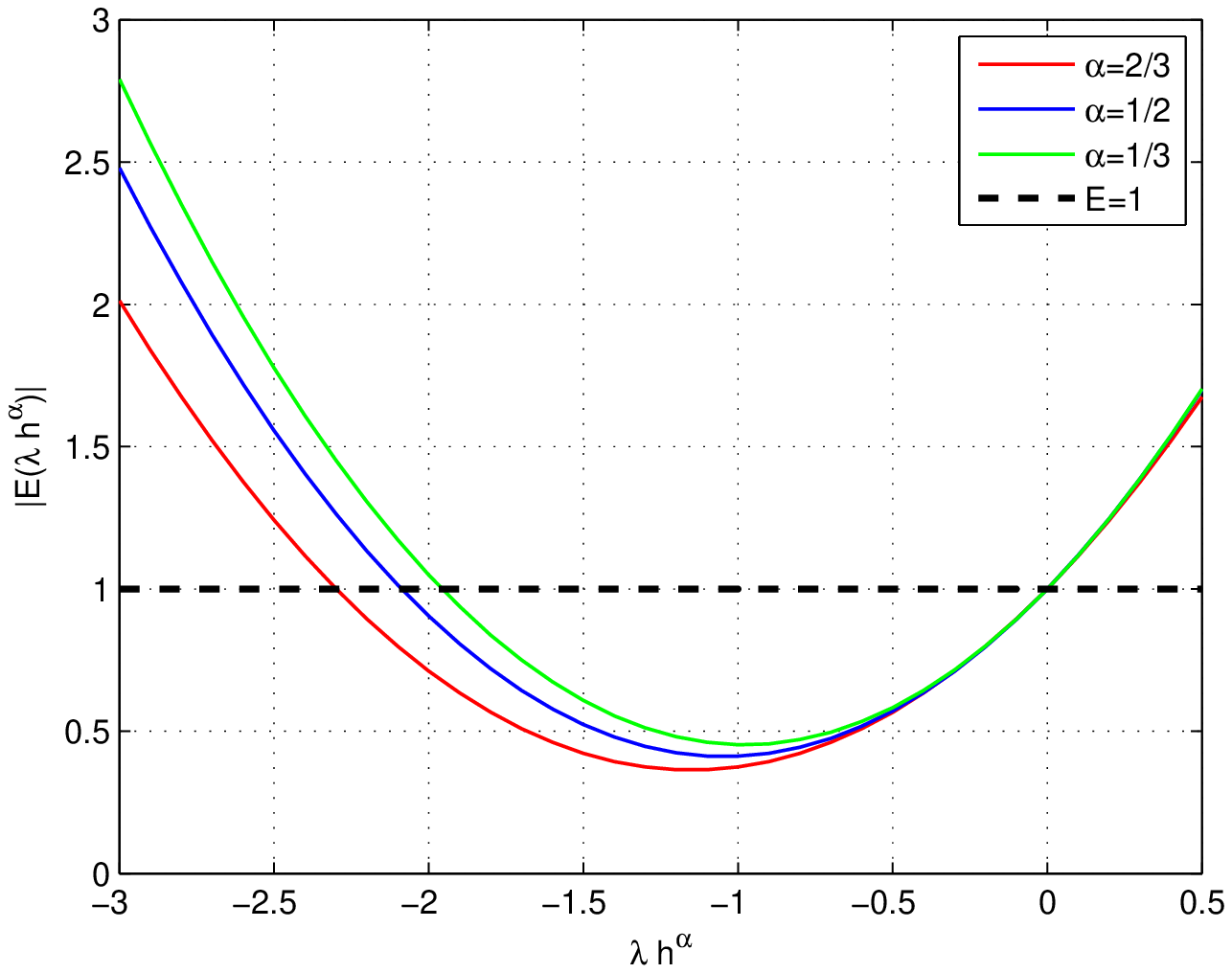}
\includegraphics[width=5.5cm]{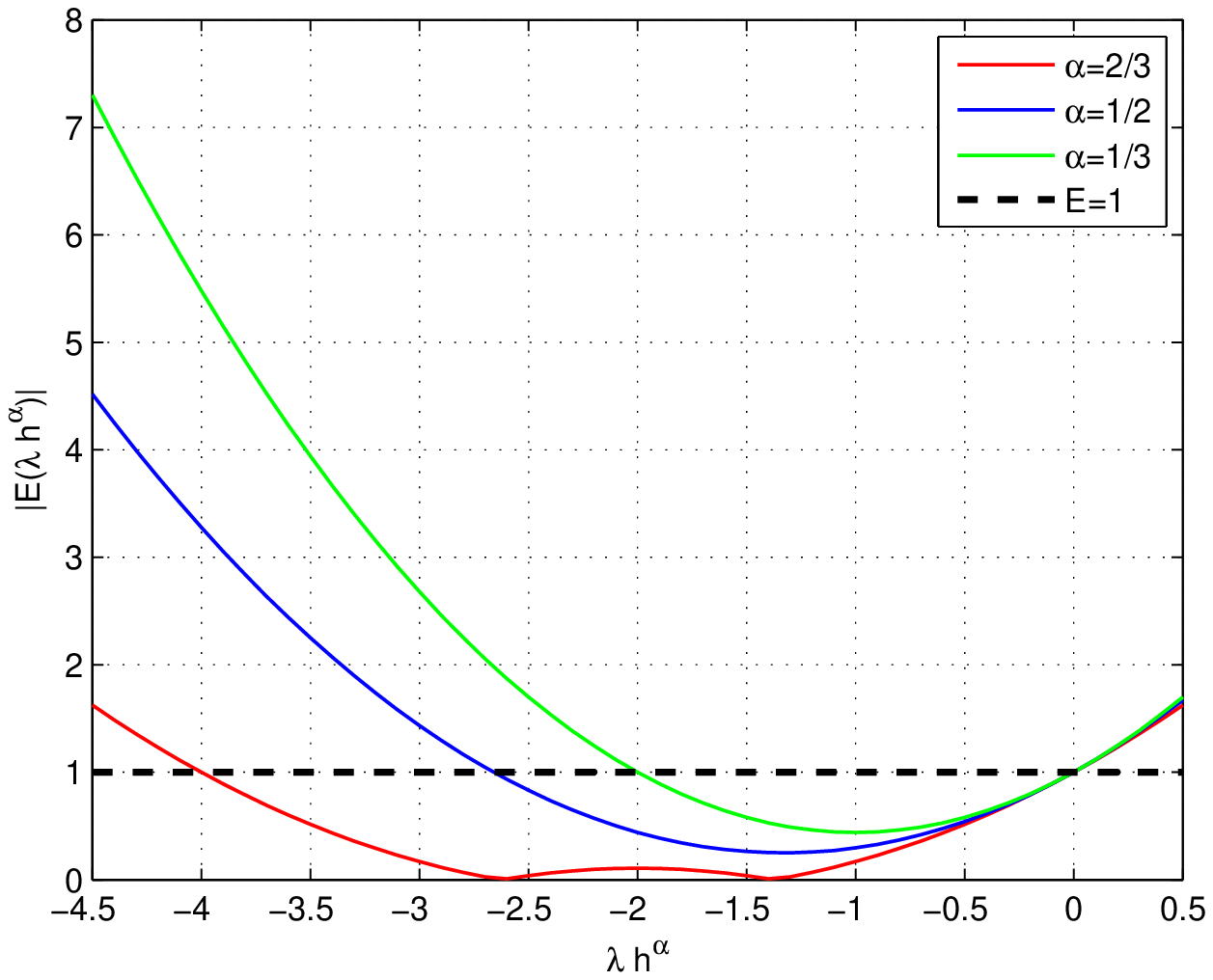}
\caption{\small{The graph of $E(\lambda h^{\alpha})$ for 2-stage
EFORK method (\ref{E10}) with
$c_2^{\alpha}=\frac{4\Gamma(\alpha+1)}{\Gamma(3\alpha+1)}$({\it
left}), and
$c_2^{\alpha}=\frac{\Gamma(\alpha+1)}{\Gamma(3\alpha+1)}$ ({\it
right}). }}
 \label{ps-2}
 \end{center}
\end{figure}

We apply the 2-stage EFORK method (\ref{E10}) to equation (\ref{s1})
and obtain
\begin{align*}
K_1&=h^{\alpha}f(t_n,y_n)= \lambda h^{\alpha}y_n,\\
K_2&=h^{\alpha}f(t_n+c_2h,y_n+a_{21}K_1)=\lambda h^{\alpha} \left( y_n+a_{21}\lambda h^{\alpha}y_n\right)\\
&=\left[ \lambda h^{\alpha}+a_{21}(\lambda h^{\alpha})^2  \right] y_n,\\
y_{n+1}&=y_n+w_1K_1+w_2K_2
=y_n+\frac{1}{2\Gamma(\alpha+1)}\left[ 2\lambda h^{\alpha} +a_{21}  (\lambda h^{\alpha})^2  \right] y_n\\
&=\left[  1+\frac{\lambda h^{\alpha}}{\Gamma(\alpha+1)} +
\frac{a_{21}(\lambda h^{\alpha})^2}{2\Gamma(\alpha+1)}  \right] y_n
=\left[  1+\frac{\lambda h^{\alpha}}{\Gamma(\alpha+1)} + \frac{c_2^{\alpha}(\lambda h^{\alpha})^2}{2(\Gamma(\alpha+1))^2}  \right] y_n\\
&=\left[  1+\frac{\lambda h^{\alpha}}{\alpha !} + \frac{c_2^{\alpha}(\lambda h^{\alpha})^2}{2(\alpha !)^2}  \right] y_n.\\
\end{align*}
Therefor, the growth factor  for 2-stage EFORK method (\ref{E10}) is
\cite{4}
\begin{equation*}\label{S1}
E(\lambda h^{\alpha})=1+\frac{\lambda h^{\alpha}}{\alpha !} +
\frac{c_2^{\alpha}(\lambda h^{\alpha})^2}{2(\alpha !)^2}\,,
\end{equation*}
So, 2-stage method (\ref{E10}) is absolutely stable if
\begin{equation*}
|1+\frac{\lambda h^{\alpha}}{\alpha !} + \frac{c_2^{\alpha}(\lambda
h^{\alpha})^2}{2(\alpha !)^2}|\leqslant 1  \,.
\end{equation*}
If $\lambda h^{\alpha}<0$, we can find the interval of absolute
stability as follows:
\begin{equation}\label{c11}
\frac{-2\alpha !}{c_2^{\alpha}}\leqslant  \lambda h^{\alpha} <0.
\end{equation}
According to (\ref{c11}), the interval of absolute stability for 2-stage EFORK method (\ref{E10}) depends on $c_2^{\alpha}$.\\
For instance, if $c_2^{\alpha}=\frac{2(\alpha !)^2}{(2\alpha)!}$\,,
then
%\begin{equation}\label{S1}
%E(\lambda h^{\alpha})=1+\frac{\lambda h^{\alpha}}{\alpha !} +
%\frac{(\lambda h^{\alpha})^2}{(2\alpha )!}\,,
%\end{equation}
%\begin{equation}
%|1+\frac{\lambda h^{\alpha}}{\alpha !} + \frac{(\lambda
%h^{\alpha})^2}{(2\alpha )!}|\leqslant 1\,,
%\end{equation}
and so the interval of absolute stability will be
\begin{equation*}
\frac{-(2\alpha) !}{\alpha !}\leqslant  \lambda h^{\alpha} <0.
\end{equation*}
Also, for
$c_2^{\alpha}=\frac{(\Gamma(2\alpha+1))^2}{\Gamma(3\alpha+1)\Gamma(\alpha+1)}$\,,
we have
%\begin{equation}
%E(\lambda h^{\alpha})=1+\frac{\lambda h^{\alpha}}{\alpha !} +
%\frac{((2\alpha)!)^2(\lambda h^{\alpha})^2}{2(3\alpha)!(\alpha
%!)^3},
%\end{equation}
%and so
\begin{equation*}
\frac{-2(\alpha !)^2(3\alpha)!}{((2\alpha)!)^2}\leqslant  \lambda
h^{\alpha} <0.
\end{equation*}
%%%%%%%
If we choose
$c_2^{\alpha}=\frac{4\Gamma(\alpha+1)}{\Gamma(3\alpha+1)}$\,, we get
%\begin{equation}
%E(\lambda h^{\alpha})=1+\frac{\lambda h^{\alpha}}{\alpha !} +
%\frac{2(\lambda h^{\alpha})^2}{(3\alpha)!(\alpha !)},
%\end{equation}
%and
\begin{equation*}
\frac{-(3\alpha)!}{2}\leqslant  \lambda h^{\alpha} <0.
\end{equation*}
%%%%%%%%%%%
or, $c_2^{\alpha}=\frac{\Gamma(\alpha+1)}{\Gamma(3\alpha+1)}$\,, we
obtain
%\begin{equation}
%E(\lambda h^{\alpha})=1+\frac{\lambda h^{\alpha}}{\alpha !} +
%\frac{(\lambda h^{\alpha})^2}{2(3\alpha)!(\alpha !)},
%\end{equation}
%and
\begin{equation*}
-2(3\alpha)!\leqslant  \lambda h^{\alpha} <0.
\end{equation*}
The graph of $E(\lambda h^{\alpha})$ for different 2-stage EFORK
methods are shown in Figures \ref{ps-1}, \ref{ps-2}. From these
Figures for $(\lambda<0)$, we can find the interval of absolute
stability for various $\alpha$.
%%%%%%%%%%%%%%%
%%%%%%%%%%%

%%%%%%%%%%%%%%%%%%%%
%%%%%%%%%%%%%%%%

\begin{figure}
\begin{center}
\includegraphics[width=6cm]{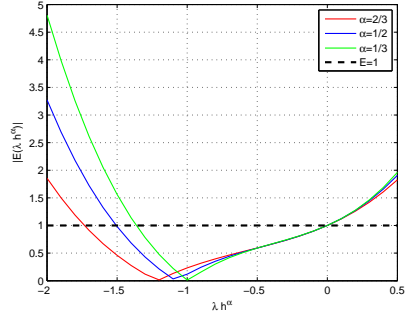}
\caption{\small{ The graph of $E(\lambda h^{\alpha})$ for 3-stage
EFORK method (\ref{E15}).}} \label{ps-3}
\end{center}
\end{figure}

Also, we apply the 3-stage EFORK method (\ref{E15}) to equation
(\ref{s1}) and get
%\begin{align*}
%K_1&=h^{\alpha}f(t_n,y_n)= \lambda h^{\alpha}y_n,\\
%K_2&=h^{\alpha}f(t_n+c_2h,y_n+a_{21}K_1)=\lambda h^{\alpha} \left( y_n+a_{21}\lambda h^{\alpha}y_n\right) \\
%&=\left[ \lambda h^{\alpha}+a_{21}(\lambda h^{\alpha})^2  \right] y_n,\\
%K_3&=h^{\alpha}f(t_n+c_3h,y_n+a_{31}K_1+a_{32}K_2)\\
%&= \left[ \lambda h^{\alpha}+a_{31}(\lambda h^{\alpha})^2
%+a_{32}(\lambda h^{\alpha})^2+a_{32}a_{21}(\lambda h^{\alpha})^3
%\right] y_n,
%\end{align*}
%with
%\begin{align*}
%y_{n+1}&=y_n+w_1K_1+w_2K_2+w_3K_3\\
%&=y_n+\left[ w_1+w_2+w_3\right] \lambda h^{\alpha} y_n+ \left[
%w_2a_{21}+w_3a_{31}+w_3a_{32} \right]   (\lambda h^{\alpha})^2y_n\\
%& \quad +\left[  w_3a_{32}a_{21}  \right] (\lambda h^{\alpha})^3y_n\,.
%\end{align*}
%By using (\ref{E17}), we have
\begin{align*}
y_{n+1}&=\left[ 1+\frac{\lambda h^{\alpha}}{\alpha !}+\frac{(\lambda
h^{\alpha})^2}{(2\alpha )!} +\frac{(\lambda h^{\alpha})^3}{(3\alpha)
!}\right] y_n\,.
\end{align*}
Thus, the growth factor for 3-stage EFORK method (\ref{E15}) is
\begin{equation*}\label{rks3}
E(\lambda h^{\alpha})= 1+\frac{\lambda h^{\alpha}}{\alpha
!}+\frac{(\lambda h^{\alpha})^2}{(2\alpha )!} +\frac{(\lambda
h^{\alpha})^3}{(3\alpha) !}\,.
\end{equation*}
The 3-stage EFORK method is absolutely stable if
\begin{equation*}
|    1+\frac{\lambda h^{\alpha}}{\alpha !}+\frac{(\lambda
h^{\alpha})^2}{(2\alpha) !} +\frac{(\lambda h^{\alpha})^3}{(3\alpha)
!}     |\leqslant 1.
\end{equation*}
The graph of $E(\lambda h^{\alpha})$ for 3-stage EFORK method
(\ref{E15}) is shown in Figure \ref{ps-3}. In this Figure, we can
see the interval of absolute stability for various $\alpha$.
\\Next, we apply the IFORK method (\ref{1i}) to equation (\ref{s1}) and get
\begin{equation*}
E(\lambda h^{\alpha})=1+\frac{\frac{1}{\alpha !}\lambda
h^{\alpha}}{1-\lambda h^{\alpha}\frac{\alpha !}{(2\alpha)!}},
\end{equation*}
with the interval of absolute stability  $(-\infty,0)$, $\lambda
<0$. In a similar manner, we can obtain the interval of absolute
stability for IFORK method (\ref{2i}). As we can see, in implicit
fractional RK methods, interval of absolute stability is very large
and they are A stable.
%%%%%%%%%%%%%%%%%%%%%%% EXAMPLES
 %%%%%%%%%%%%%%%%%%%%%%%%%%%%%%%%%%%%%%%%%%%%%%%%%%%%%%%%%%%
 \section{Numerical examples}\label{exam}
 In order to demonstrate the effectiveness and order of accuracy of the proposed methods in sections \ref{se1}-\ref{se2}, two examples are considered.\\
\\
\\
\\
 $\bf{Example}$ $\bf{1:}$ Let us consider fractional differential
 equation
 \begin{eqnarray*}
&& {}_{0}^cD_t^{\alpha} y(t)=-y(t)+\frac{t^{4-\alpha}}{\Gamma(5-\alpha)},\quad t\in[ 0,T],\nonumber\\
&& y(0)=0,
\end{eqnarray*}
 where, the exact solution of the equation is $y(t)=t^4E_{\alpha,5}(-t^{\alpha})$. \\
 For different values of $h$, $\alpha$ and $T$, the computed solutions are compared with the exact solution.
 We have reported the absolute error in time $T$ as:
$$ E(h,T)=|y(t_{N^m})-y_{N^m}|. $$
Also, we calculated the computational orders of the presented method
according to the following relation:
$$Log_2\frac{E(h,T)}{E(h/2,T)}.$$

The computed solutions by 2-stage EFORK method (\ref{E10}), 3-stage EFORK method (\ref{E15}) and IFORK methods (\ref{1i})
and (\ref{2i}) are reported in Tables \ref{T1}-\ref{IT3}. From the Tables \ref{T1}-\ref{IT3}, we can conclude that
the computed orders of  truncation errors is in a good agreement with the obtained results  of sections \ref{se1}-\ref{se2} .\\

 %%%%%%%%%%%%%%%%%%%%%%%%%%%%%%%%%%%%%
 \begin{table}[!h]
\begin{center}
\caption{\small{ 2-stage method (\ref{E10}) for $T=1$, $m=3$ in
Example 1.}}
\begin{tabular}{c c c c}
\hline$\alpha$  &\qquad$h$  &\qquad $E(h,T)$   &\qquad $Log_2\frac{E(h,T)}{E(h/2,T)} $                    \\
\hline 1/3 &\qquad 1/40 &\qquad $1.09027e{-2}$  &\qquad $ 1.1320$            \\
 &\qquad 1/80         &\qquad $4.97465e{-3}$       &\qquad   $1.0013$    \\
&\qquad 1/160        &\qquad  $2.48509e{-3}$         &\qquad    $0.9136$\\
 &\qquad 1/320           &\qquad$1.31920e{-3}$          &\qquad    $ 0.8533$  \\
 &\qquad 1/640             &\qquad$7.30171e{-4}$         &\qquad    $*$ \\
\hline
\end{tabular}\label{T1}
\end{center}
\end{table}
  %%%%%%%%%%%%%%%%%%%%%%%%%%%%%%%%%%%%%%%%%%
    \begin{table}[!h]
\begin{center}
\caption{\small{ 2-stage method (\ref{E10}) for $T=1$, $m=2$ in
Example 1.}}
\begin{tabular}{c c c c}
\hline$\alpha$  &\qquad$h$  &\qquad $E(h,T)$    &\qquad $Log_2\frac{E(h,T)}{E(h/2,T)} $                    \\
\hline 1/2 &\qquad 1/40 &\qquad$2.05503e{-3}$ &\qquad $1.2248$            \\
 &\qquad 1/80          &\qquad$8.79256e{-4}$     &\qquad   $1.1621$    \\
&\qquad 1/160       &\qquad$3.92907e{-4}$       &\qquad    $   1.1171$\\
 &\qquad 1/320        &\qquad$1.81137e{-4}$          &\qquad    $ 1.0845$  \\
 &\qquad 1/640             &\qquad$8.54183e{-5}$         &\qquad    $*$ \\
\hline
\end{tabular}\label{T2}
\end{center}
\end{table}
%%%%%%%%%%%%%%%%%%%%%%%%%%%%
 %%%%%%%%%%%%%%%%%%%%%%%%%%%%%%%%%%%%%%%%%%%
   \begin{table}[!h]
\begin{center}
\caption{\small{ 3-stage method (\ref{E15}) for $T=1$, $m=4$ in
Example 1.}}
\begin{tabular}{c c c c}
\hline$\alpha$  &\qquad$h$  &\qquad $E(h,T)$    &\qquad $Log_2\frac{E(h,T)}{E(h/2,T)} $                    \\
\hline 1/4 &\qquad 1/40  &\qquad $9.94252e{-4}$ &\qquad $0.8437$            \\
 &\qquad 1/80         &\qquad $5.54011e{-4}$         &\qquad   $ 0.8214$    \\
&\qquad 1/160         &\qquad  $3.13499e{-4}$           &\qquad    $0.8064$\\
 &\qquad 1/320           &\qquad$1.79258e{-4}$          &\qquad    $ 0.7958$  \\
 &\qquad 1/640             &\qquad$1.03255e{-4}$         &\qquad    $*$ \\
\hline
\end{tabular}\label{TT3}
\end{center}
\end{table}
%%%%%%%%%%%%%%%%%%%%%%%%%%%%%%%%%%%%%%%%
 %%%%%%%%%%%%%%%%%%%%%%%%%%%%%%%
  \begin{table}[!h]
\begin{center}
\caption{\small{ 3-stage method (\ref{E15}) for $T=1$, $m=2$ in
Example 1.}}
\begin{tabular}{c c c c}
\hline$\alpha$  &\qquad$h$  &\qquad $E(h,T)$    &\qquad $Log_2\frac{E(h,T)}{E(h/2,T)} $                    \\
\hline 1/2 &\qquad 1/40  &\qquad $7.45694e{-5}$ &\qquad $1.5942$            \\
 &\qquad 1/80         &\qquad $2.46986e{-5}$         &\qquad   $1.5789$    \\
&\qquad 1/160         &\qquad  $8.26771e{-6}$           &\qquad    $ 1.5625$\\
 &\qquad 1/320           &\qquad$2.79911e{-6}$          &\qquad    $ 1.5478$  \\
 &\qquad 1/640             &\qquad$9.57367e{-7}$         &\qquad    $*$ \\
\hline
\end{tabular}\label{T3}
\end{center}
\end{table}
 %%%%%%%%%%%%%%%%%%%%%%%%%%%%%%%%%%%
 %%%%%%%%%%%%%%%%%%%%%%%%%%%%%%%%%%%%%
%%%%%%%%%%%%%%%%%%%%%%%%%%
 \begin{table}[!h]
\begin{center}
\caption{\small{IFORK methods for $T=1$, $m=2$ in Example 1.}}
\begin{tabular}{c c c c c c}
\hline$\alpha$  &\hspace{-.071 cm} $h$  &  IFORK (\ref{1i})    & $Log_2\frac{E(h,T)}{E(h/2,T)} $ & IFORK (\ref{2i}) & $Log_2\frac{E(h,T)}{E(h/2,T)} $              \\
\hline 1/2  &  1/40   \hspace{-.071 cm}      & $1.86448e{-4}$    & $1.0832$      &  $3.56080e{-4}$   &   $1.7307$ \\
                  & 1/80  \hspace{-.071 cm}        & $8.79978e{-5}$    &  $ 1.0492$    &   $1.07291e{-4}$   &   $1.6666$               \\
                  & 1/160   \hspace{-.071 cm}     &$4.25221e{-5}$      &  $1.0291$     &   $3.37953e{-5}$    &   $1.6191$          \\
                  & 1/320   \hspace{-.071 cm}     &  $2.08361e{-5}$     & $ 1.0174$    &  $1.10016e{-5}$ &      $1.5846$           \\
                  & 1/640   \hspace{-.071 cm}      & $1.02928e{-5}$     & $*$              &  $3.66806e{-6}$           &     $*$             \\
\hline
\end{tabular}\label{IT2}
\end{center}
\end{table}
 %%%%%%%%%%%%%%%%%%%%%%%%%%%%%%%%%%%%%%%%%%%%%%%%%
   \begin{table}[!h]
\begin{center}
\caption{\small{IFORK methods for $T=1$, $m=3$ in Example 1.}}
\begin{tabular}{c c c c c c}
\hline$\alpha$  &\hspace{-.071 cm} $h$  &  IFORK (\ref{1i})    & $Log_2\frac{E(h,T)}{E(h/2,T)} $ & IFORK (\ref{2i}) & $Log_2\frac{E(h,T)}{E(h/2,T)} $              \\
\hline 1/3 & 1/40   \hspace{-.071 cm}     & $3.04275e{-4}$ &  $ 0.9083$    &  $5.12337e{-3}$   &    $1.4661$ \\
           & 1/80   \hspace{-.071 cm}     & $1.62123e{-4}$     &   $ 0.8746$    &   $1.85445e{-3}$   &    $1.3856$               \\
           & 1/160   \hspace{-.071 cm}    & $8.84207e{-5}$       &    $ 0.8482$ &    $7.09751e{-4}$    &   $1.2940$          \\
           & 1/320   \hspace{-.071 cm}     & $4.91160e{-5}$          &    $0.8255$   &  $2.89458e{-4}$ &       $1.2106$           \\
           & 1/640   \hspace{-.071 cm}     & $2.77152e{-5}$         &    $*$  &     $1.25070e{-4}$           &     $*$             \\
\hline
\end{tabular}\label{IT3}
\end{center}
\end{table}

\begin{figure}
\begin{center}
\includegraphics[width=5cm]{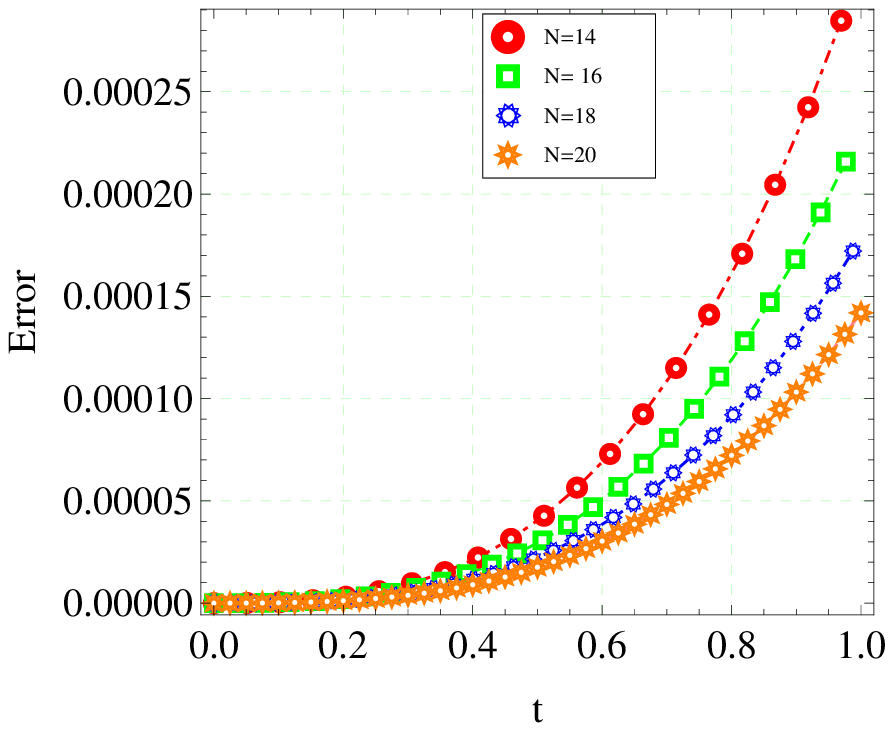}
\includegraphics[width=5cm]{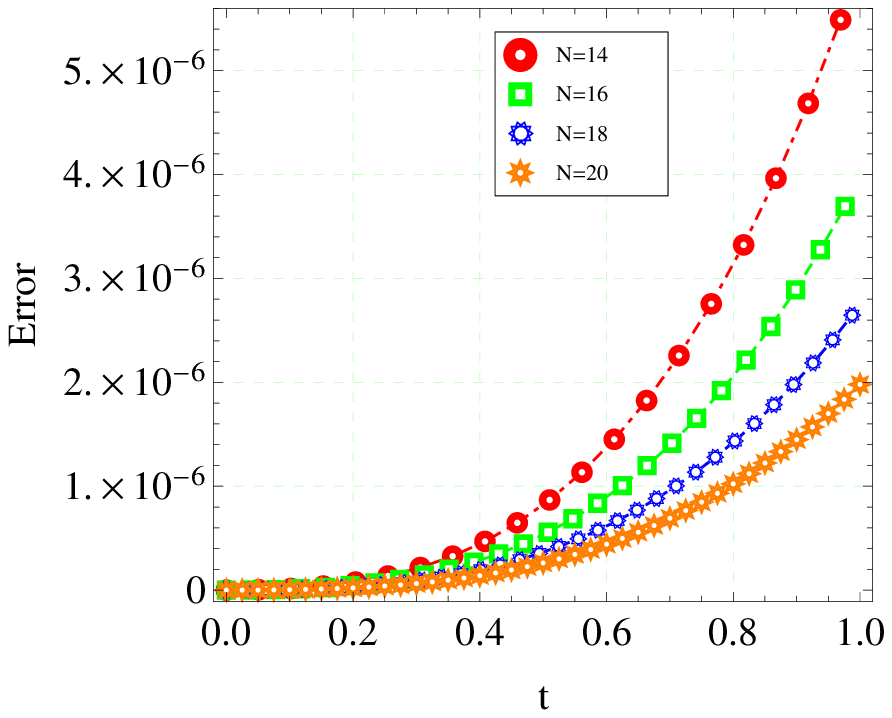}
\caption{\small{ The error curves of the 2-stage EFORK method
(\ref{E10}) ({\it left}), and the error curves of the 3-stage EFORK
method (\ref{E15})  ({\it right}) for Example 1. }}
 \label{er-1}
 \end{center}
\end{figure}

Fig \ref{er-1}, illustrates the error curves of the 2-stage EFORK
method (\ref{E10}) and  the 3-stage EFORK method (\ref{E15}) at
$T=1,$ with $\alpha=1/2$, $m=2$ and different values of $N$.\\\\
%%%%%%%%%%%%%%%%%%%%%%%%%%%%%%%%%%%%%%
\vspace{0 cm}
\\
\\
\\
 %%%%%%%%%%%%%%%%%%%%%%%%%%%%%%%%%%%%%%%%%%%%%%%%%%%%%%
\\$\bf{Example}$ $\bf{2:}$ Let us consider the following fractional differential equation from \cite{9}:
 \begin{eqnarray*}
&& {}_{0}^cD_t^{\alpha} y(t)=\frac{2}{\Gamma(3-\alpha)}t^{2-\alpha}-\frac{1}{\Gamma(2-\alpha)}t^{1-\alpha}-y(t)+t^2-t,\quad t\in[ 0,T],\nonumber\\
&& y(0)=0,
\end{eqnarray*}
 where the exact solution of the problem is $y(t)=t^2-t$. \\
 For different values of $h$, $\alpha$ and $T$, the obtained results by 2-stage EFORK method (\ref{E10}), 3-stage EFORK method (\ref{E15}).
 IFORK methods (\ref{1i}) and (\ref{2i}) are shown in Tables \ref{T5}-\ref{IT6}.

  Also, from the Tables \ref{T5}-\ref{IT6}, we can see that the computed orders  is consistent with the given results of  sections \ref{se1}-\ref{se2}.\\

\begin{figure}
\begin{center}
\includegraphics[width=5cm]{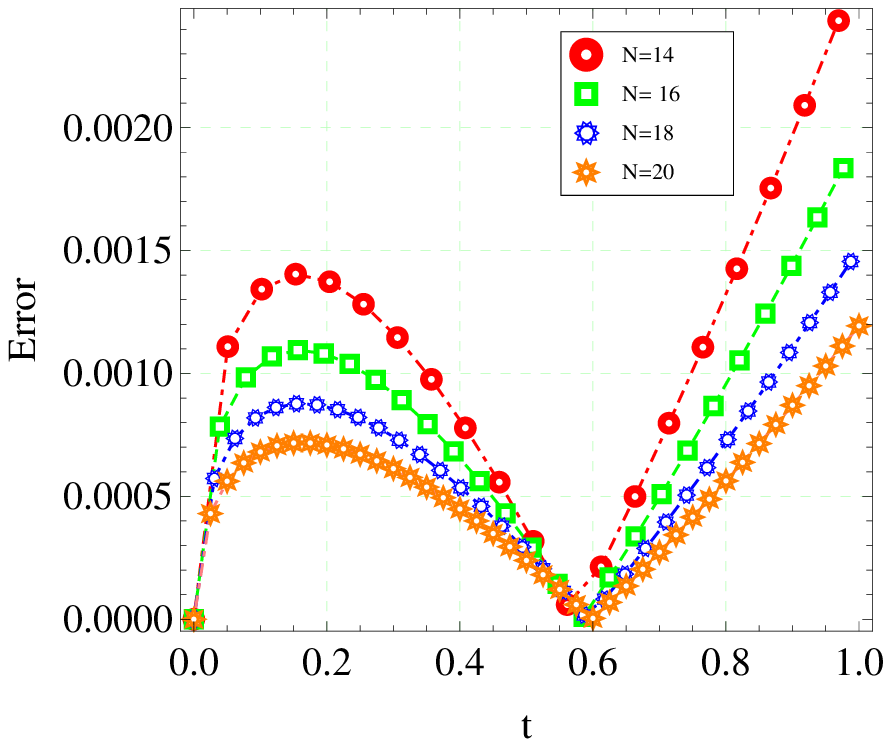}
\includegraphics[width=5cm]{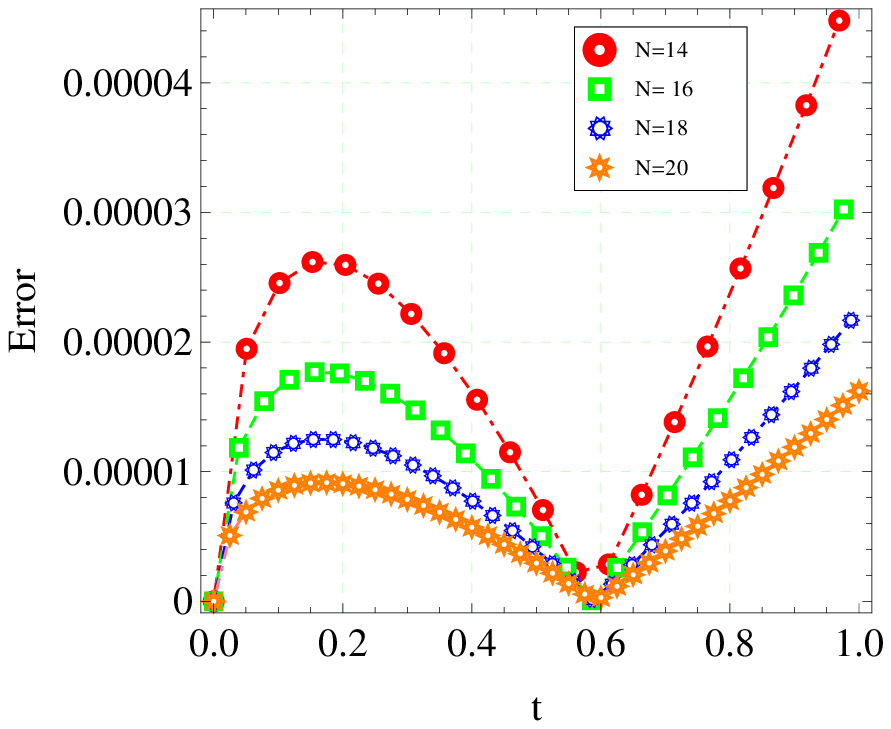}
\caption{\small{ The error curves of the 2-stage EFORK method
(\ref{E10}) ({\it left}), and the 3-stage EFORK method (\ref{E15})
({\it right}) in Example 2. }}
 \label{er-2}
 \end{center}
\end{figure}

\begin{figure}
\begin{center}
\includegraphics[width=5cm]{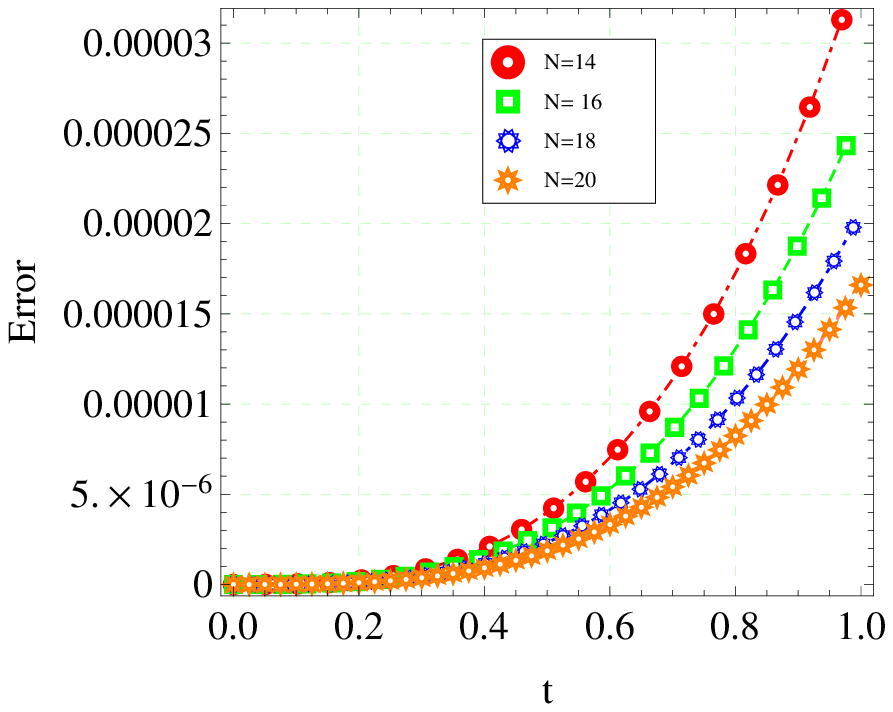}
\includegraphics[width=5cm]{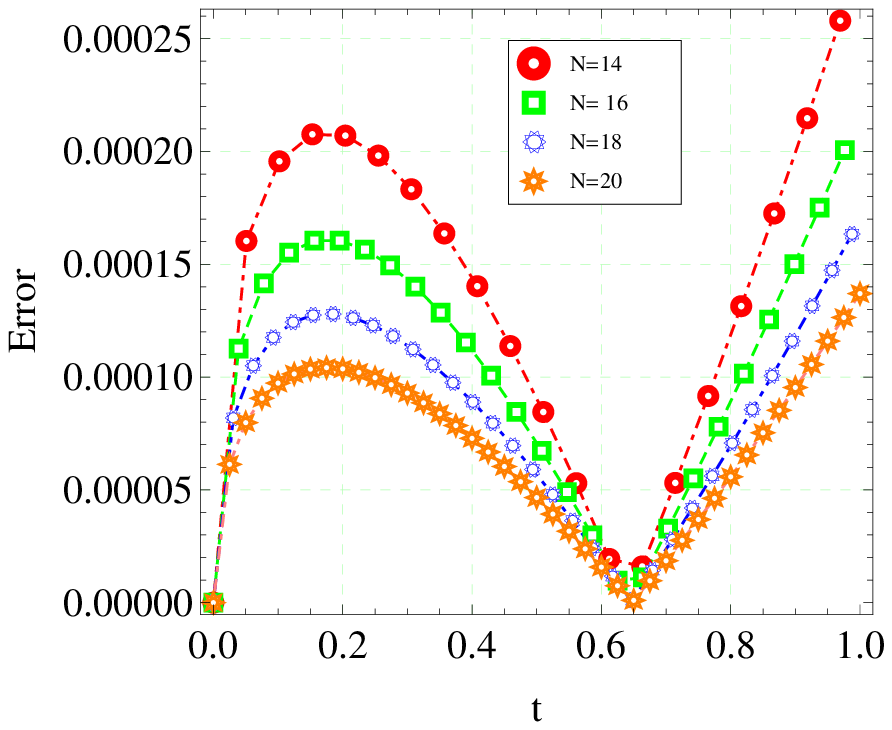}
\caption{\small{ The error curves of the IFORK method (\ref{1i}) in
Example 1 ({\it left}), and  Example 2 ({\it right}). }}
 \label{er-3}
 \end{center}
\end{figure}

  Fig \ref{er-2}, illustrates the numerical results of the 2-stage EFORK method (\ref{E10}) and the 3-stage EFORK method (\ref{E15})
  at $T=1$ for $\alpha=1/2$, $m=2$ and different values of $N$.
  Also, Fig \ref{er-3} illustrates the numerical results of the IFORK method (\ref{1i}) for Example 1 and  Example 2  at $T=1$
  for $\alpha=1/2$, $m=2$ and different values of $N$.
 %%%%%%%%%%%%%%%%%%%%%%%%%%%
  \begin{table}[!h]
\begin{center}
\caption{\small{ 2-stage method (\ref{E10}) for $T=1$, $m=3$ in
Example 2.}}
\begin{tabular}{c c c c}
\hline$\alpha$  &\qquad$h$  &\qquad $E(h,T)$    &\qquad $Log_2\frac{E(h,T)}{E(h/2,T)} $                    \\
\hline 1/3 &\qquad 1/40 &\qquad $1.00356e{-1}$   &\qquad $ 1.1075$            \\
 &\qquad 1/80         &\qquad $4.65748e{-2}$        &\qquad   $ 0.9928$    \\
&\qquad 1/160        &\qquad   $2.34046e{-2}$       &\qquad    $   0.9107$\\
 &\qquad 1/320           &\qquad$1.24493e{-2}$            &\qquad    $   0.8523$  \\
 &\qquad 1/640             &\qquad$6.89556e{-3}$          &\qquad    $*$ \\
\hline
\end{tabular}\label{T5}
\end{center}
\end{table}
 %%%%%%%%%%%%%%%%%%%%%%%%%%%%%%%%
 \begin{table}[!h]
\begin{center}
\caption{\small{ 2-stage method (\ref{E10}) for $T=1$, $m=2$ in
Example 2.}}
\begin{tabular}{c c c c}
\hline$\alpha$  &\qquad$h$  &\qquad $E(h,T)$   &\qquad $Log_2\frac{E(h,T)}{E(h/2,T)} $                    \\
\hline 1/2 &\qquad 1/40 &\qquad $1.77152e{-2}$   &\qquad $  1.2351$            \\
 &\qquad 1/80         &\qquad $7.52581e{-3}$         &\qquad   $  1.1738$    \\
&\qquad 1/160        &\qquad   $3.33574e{-3}$     &\qquad    $   1.1275$\\
 &\qquad 1/320           &\qquad$1.52680e{-3}$           &\qquad    $  1.0928$  \\
 &\qquad 1/640             &\qquad$7.15859e{-4}$          &\qquad    $*$ \\
\hline
\end{tabular}\label{T6}
\end{center}
\end{table}
 %%%%%%%%%%%%%%%%%%%%%%%%%%%%
 \begin{table}[!h]
\begin{center}
\caption{\small{ 3-stage method (\ref{E15}) for $T=1$, $m=4$ in
Example 2.}}
\begin{tabular}{c c c c}
\hline$\alpha$  &\qquad$h$  &\qquad $E(h,T)$     &\qquad $Log_2\frac{E(h,T)}{E(h/2,T)} $                        \\
\hline 1/4 &\qquad 1/40  &\qquad $9.90939e{-3}$      &\qquad     $  0.8383$  \\
 &\qquad 1/80         &\qquad $5.54249e{-3}$         &\qquad      $0.8182$ \\
&\qquad 1/160         &\qquad  $3.14342e{-3}$        &\qquad     $  0.8047$  \\
 &\qquad 1/320           &\qquad$1.79955e{-3}$      &\qquad      $  0.7950$    \\
 &\qquad 1/640             &\qquad$1.03718e{-3}$     &\qquad $*$         \\
\hline
\end{tabular}\label{TT7}
\end{center}
\end{table}
%%%%%%%%%%%%%%%%%%%%%%%%%%%%%%
  %%%%%%%%%%%%%%%%%%%%%%%%%%%%%
 \begin{table}[!h]
\begin{center}
\caption{\small{ 3-stage method (\ref{E15}) for $T=1$, $m=2$ in
Example 2.}}
\begin{tabular}{c c c c}
\hline$\alpha$  &\qquad$h$  &\qquad $E(h,T)$     &\qquad $Log_2\frac{E(h,T)}{E(h/2,T)} $                        \\
\hline 1/2 &\qquad 1/40  &\qquad $5.79341e{-4}$      &\qquad     $  1.5592$  \\
 &\qquad 1/80         &\qquad $1.96590e{-4}$         &\qquad      $1.5566$ \\
&\qquad 1/160         &\qquad  $6.68302e{-5}$        &\qquad     $  1.5475$  \\
 &\qquad 1/320           &\qquad$2.28624e{-5}$      &\qquad      $ 1.5374$    \\
 &\qquad 1/640             &\qquad$7.87606e{-6}$     &\qquad $*$         \\
\hline
\end{tabular}\label{T7}
\end{center}
\end{table}
%%%%%%%%%%%%%%%%%%%%%%%%%%%%%%%%
  \begin{table}[!h]
\begin{center}
\caption{\small{IFORK methods for $T=1$, $m=2$ in Example 2.}}
\begin{tabular}{c c c c c c}
\hline$\alpha$  &\hspace{-.071 cm} $h$  &  IFORK (\ref{1i})    & $Log_2\frac{E(h,T)}{E(h/2,T)} $ & IFORK (\ref{2i}) & $Log_2\frac{E(h,T)}{E(h/2,T)} $              \\
\hline   1/2 & 1/40  \hspace{-.071 cm}         &$1.52888e{-3}$        & $ 1.0777$     &  $2.99223e{-3}$    &    $1.6608$ \\
                 &  1/80   \hspace{-.071 cm}       & $7.24362e{-4}$       &   $1.0464$    &   $9.46328e{-4}$   &    $1.6244$               \\
                 &1/160   \hspace{-.071 cm}       &$3.50728e{-4}$        &    $1.0279$   &    $3.06932e{-4}$   &  $1.5943$          \\
                & 1/320     \hspace{-.071 cm}     &$1.72002e{-4}$         &   $1.0171$   &  $1.01649e{-4}$      &  $1.5703$           \\
                 & 1/640     \hspace{-.071 cm}    &$8.49858e{-5}$         &    $*$          &    $3.42297e{-5}$     &     $*$             \\
\hline
\end{tabular}\label{IT5}
\end{center}
\end{table}
%%%%%%%%%%%%%%%%%%%%%%%%%%%%%%%%%%%%%%
 \begin{table}[!h]
\begin{center}
\caption{\small{IFORK methods for $T=1$, $m=3$ in Example 2.}}
\begin{tabular}{c c c c c c}
\hline$\alpha$  &\hspace{-.071 cm} $h$  &  IFORK (\ref{1i})    & $Log_2\frac{E(h,T)}{E(h/2,T)} $ & IFORK (\ref{2i}) & $Log_2\frac{E(h,T)}{E(h/2,T)} $              \\
\hline  1/3 & 1/40    \hspace{-.071 cm}               &$2.84493e{-3}$      & $0.9005$    &  $4.77620e{-2}$   &    $1.4860$ \\
                & 1/80     \hspace{-.071 cm}              &$1.52403e{-3}$     &   $0.8707$    &   $1.70509e{-2}$   &    $1.3835$               \\
               & 1/160     \hspace{-.071 cm}              &$8.33483e{-4}$     &    $0.8463$  &    $6.53520e{-3}$    &   $1.2859$          \\
              & 1/320     \hspace{-.071 cm}           &$4.63592e{-4}$       &    $0.8246$   &  $2.68020e{-3}$         &      $1.2023$           \\
              & 1/640        \hspace{-.071 cm}         &$2.61754e{-4}$         &    $*$            &     $1.16478e{-3}$       &     $*$             \\
\hline
\end{tabular}\label{IT6}
\end{center}
\end{table}
%%%%%%%%%%%%%%%%%%%%%%%%%%%%%%%
 %%%%%%%%%%%%%%%%%%%%%%%%%%%%%%%%%%%%%%%%%%%%
 \begin{table}[!h]
\begin{center}
\caption{\small{$E(h,T)$ for $\alpha=1/2$, $m=2$ and different
values of $T$.}}
\begin{tabular}{ c c c c c}
\hline$T$ &\quad 2-stage-Exam.1 &\quad 3-stage-Exam.1 &\quad 2-stage-Exam.2 &\quad 3-stage-Exam.2 \\
\hline 0.5 &~~ $4.81351e{-6}$ &\quad$4.20218e{-8}$ &\quad $9.29068e{-5}$ &\quad$5.57177e{-7}$ \\
1.0 &\quad $8.54183e{-5}$ &\quad$9.57367e{-7}$ &\quad $3.00335e{-3}$ &\quad$7.87606e{-6}$ \\
1.5 &\quad $4.51426e{-4}$ &\quad$6.00243e{-6}$ &\quad $2.78127e{-3}$ &\quad$3.63915e{-5}$ \\
2 &\quad$1.45778e{-3}$ &\quad$2.20455e{-5}$ &\quad$6.26033e{-3}$ &\quad$9.38735e{-5}$ \\
3 &\quad$7.50603e{-3}$ &\quad$1.37070e{-4}$ &\quad$1.78449e{-2}$ &\quad$3.26003e{-4}$ \\
\hline
\end{tabular}\label{T4}
\end{center}
\end{table}
%%%%%%%%%%%%%%%%%%%%%%%%%%%%%%%%%%%%%%%
 \begin{table}[!h]
\begin{center}
\caption{\small{ Optimal 2-stage method (\ref{E10}) for $T=1$,
$m=2$. }}
\begin{tabular}{l  c  c c}
\hline$\alpha$   & \hspace{0 cm} $h$                & $E(h,T)$, Example 1  & $E(h,T)$,  Example 2  \\
\hline 1/2       &\hspace{0 cm} 1/40 & $7.35533e{-4}$   & $6.12299e{-3}$       \\
                 & \hspace{0 cm}1/80        & $3.55401e{-4}$      & $2.94885e{-3}$         \\
                 & \hspace{0 cm}1/160      &  $1.72778e{-4}$      & $1.43060e{-3}$      \\
                 & \hspace{0 cm}1/320      & $8.45336e{-5}$       & $6.99024e{-4}$     \\
                 & \hspace{0 cm}1/640      & $4.15855e{-5}$       & $3.43589e{-4}$       \\
\hline
\end{tabular}\label{T9}
\end{center}
\end{table}
%%%%%%%%%%%%%%%%%%%%%%%%%%%%%%

%%%%%%%%%%%%%%%%%%%%%%%%%%%%%%%%%%%%%%
From these Tables, we conclude that the computational order for 2
and 3 stages EFORK methods are ${2\alpha}$ and ${3\alpha}$,
respectively. As expected and seen in Tables and Figs, the 3-stage
EFORK method in comparison with the 2-stage EFORK method, provides
better results. In Table \ref{T4}, the numerical results for
different values of $T$ are shown.
%%%%%%%%%%%%%%%%%
%%%%%%%%%%%%%%%%%%%%%%%
Also, the numerical results for optimal case
$c_2^{\alpha}=\frac{(\Gamma(2\alpha+1))^2}{\Gamma(3\alpha+1)\Gamma(\alpha+1)}$,
in 2-stage EFORK method are shown in Table \ref{T9} with
$\alpha=1/2$ and different values of $h$.\\
\\
%%%%%%%%%%%%%%%%%%%%%%%%%%%%%%%Conclusion
\section{Conclusions}
This paper introduces new efficient FORK methods for FDEs based on Caputo generalized Taylor formulas. The proposed methods were examined for consistency, convergence, and stability. The interval of absolute stability of FORK methods has been determined, and implicit fractional order RK methods were shown to be A %Please check if this is correct
stable. Some examples were provided to demonstrate the effectiveness of these numerical schemes. We can obtain these results for Riemann--Liouville and Gronwald--Letnikov fractional derivatives accordingly.\\
{Recently, a new concept of differentiation called fractal and fractional differentiation was suggested and numerically examined by many researchers \cite{461,462}, where the differential operator has two orders: the first is fractional order and the second is the fractal dimension. These differential (integral) operators have not been studied intensively yet. In future work, we will extend the presented method for fractional differential equations with fractal--fractional derivatives.}\\\\

\end{document}